\documentclass[10pt]{article}
\usepackage[utf8]{inputenc}
\usepackage[english]{babel}
\usepackage{amsmath,amsthm,amsfonts}
\usepackage{graphicx}
\usepackage{bm}
\usepackage{tikz}
\usepackage{multirow}
\usepackage{color}

\graphicspath{ {./figs/} }

\title{Nonlocal multicontinuum (NLMC) upscaling of mixed dimensional coupled flow problem for embedded and discrete fracture models}

\author{
Maria Vasilyeva \thanks{Institute for Scientific Computation, Texas A\&M University, College Station, TX 77843-3368 \& Department of Computational Technologies, North-Eastern Federal University, Yakutsk, Republic of Sakha (Yakutia), Russia, 677980. Email: {\tt vasilyevadotmdotv@gmail.com}.}
\and
Eric T. Chung \thanks{Department of Mathematics,
The Chinese University of Hong Kong (CUHK), Hong Kong SAR. Email: {\tt tschung@math.cuhk.edu.hk}.}
\and
\and
Wing Tat Leung
\thanks{Institute of Computational Engineering and Sciences, University of Texas at Austin, Austin, USA}
\and
Valentin Alekseev 
\thanks{Multiscale model reduction laboratory, North-Eastern Federal University, Yakutsk, Republic of Sakha (Yakutia), Russia, 677980.}
}

\begin{document}

\maketitle

\begin{abstract}
In this work, we present an upscaled model for mixed dimensional coupled flow problem in fractured porous media. We consider both embedded and discrete fracture models (EFM and DFM) as fine scale models which contain coupled system of equations. For fine grid discretization, we use a conservative finite-volume approximation. We construct an upscaled model using the non-local multicontinuum (NLMC) method for the coupled system. The proposed upscaled model is based on a set of simplified multiscale basis functions for the auxiliary space and a constraint energy minimization principle for the construction of multiscale basis functions. Using the constructed NLMC-multiscale basis functions, we obtain an accurate coarse grid upscaled model. We present numerical results for both fine-grid models and upscaled coarse-grid models using our NLMC method.
We consider model problems with (1) discrete fracture fine grid model with low and high permeable fractures; (2) embedded fine grid model for two types of geometries with differnet fracture networks and (3) embedded fracture fine grid model with heterogeneous permeability. The simulations using the upscaled model provide very accurate solutions with significant reduction in the dimension of the problem.
\end{abstract}

\section*{Introduction}

Mathematical simulation of the flow processes in fractured porous media plays an important role in reservoir simulation, nuclear waste disposal, $CO_2$ sequestration, unconventional gas production and geothermal energy production.
Fracture networks usually have complex geometries, multiple scales and very small thickness compared to typical reservoir sizes.
Due to high permeability, fractures have a significant impact on the flow processes. 
A common approach to model fracture media is to consider the discrete fractures as lower-dimensional objects \cite{martin2005modeling, d2012mixed, formaggia2014reduced, Quarteroni2008coupling, schwenck2015dimensionally}.
This results in a coupled mixed dimensional mathematical models, where we have $d$ - dimensional equation for flow in porous matrix and $(d-1)$ - dimensional equation for fracture networks. The cross-flow equilibrium between the flow in fracture and matrix is described by some specific transfer terms.

Various numerical approaches to model fractured porous media have been developed and can be classified by the types of meshing techniques used for simulations.
One approach, called discrete fracture model (DFM) is associated with the conforming discretization or explicit meshing of the fracture geometry. 
In DFM, we consider flow processes in both the matrix and the fracture media, where the fractures are located at the interfaces between matrix cells \cite{hoteit2008efficient, karimi2003efficient, karimi2001numerical, garipov2016discrete}.
This requires a conforming meshing of the fractures, which can lead to large computational demands since a large number of unknowns arise when the problem is approximated. Nevertheless, DFM is shown to be an accurate tool to describe the flow characteristics of the models with large-scale fractures.
In another approach, called the embedded fracture model (EFM) \cite{hkj12, ctene2016algebraic, tene2016multiscale} the fractures are not resolved by grid but are considered as an overlaying continua. In EFM, matrix and fracture are viewed as two porosity types co-existing at the same spatial location, thus simple structured meshes can be used for the domain discretization. The transfer term between matrix and fracture appears as an additional source or sink and is assumed to exist in entire simulation domain. The concept of this approach can be classified in the class of dual-continuum or multi-continuum models \cite{barenblatt1960basic, warren1963behavior, douglas1990dual, ginting2011application}. 

In this work, we consider both embedded and discrete fracture models (EFM and DFM) for fine-scale model. Mathematical models for both approaches are described by the coupled mixed dimensional system.
Finite volume methods are widely used discretization techniques for simulation of flow problems. 
For fine grid simulations, we employ the cell-centered finite-volume method with two-point flux approximation (TPFA) \cite{hkj12, ctene2016algebraic, bosma2017multiscale, ctene2017projection, tene2016multiscale}.
In the DFM approach, we impose Robin type conditions on the internal boundaries that represent fractures. This allows us to couple the subdomains by simply discretizing the flux over faces of each internal boundary.
In the EFM approach, a coupling between fracture and matrix subdomains is described by some transfer term.

Due to the scale disparity, fine grid simulation of the processes in fractured porous media can be very expensive if one needs to capture various scales of flow features at once.
To reduce the dimension of the fine scale system directly using finite volume approximation of the problem with EFM and DFM approaches, multiscale methods or upscaling techniques are proposed \cite{houwu97, eh09, weinan2007heterogeneous, lunati2006multiscale, jenny2005adaptive}.
In our previous work, the multiscale model reduction techniques based on the Generalized multiscale finite element method (GMsFEM) for flow in fractured porous media are presented \cite{akkutlu2015multiscale, chung2017coupling, efendiev2015hierarchical, akkutlu2018multiscale}. The general idea of GMsFEM is to first solve some local problems to get snapshot spaces, then design suitable spectral problems to obtain important modes which can be used to construct multiscale basis \cite{EGG_MultiscaleMOR, egh12, chung2016adaptive, CELV2015}. The resulting multiscale space contains basis functions that take into account the microscale heterogeneities, and the multiscale scale solution found in this space provide an accurate approximation. Recently, the authors in \cite{chung2017constraint} proposed a new GMsFEM method with constraint energy minimization (CEM-GMsFEM). In CEM-GMsFEM, one constructs multiscale basis functions which can capture long channelized effects and can be localized in an oversampling domain. The construction of the multiscale space starts with an auxiliary space, which consists of eigenfunctions of local spectral problems. Using the auxiliary space, one can obtain the required multiscale spaces by solving a constraint energy minimization problem. Using the multiscale basis functions, we recently presented a non-local multi-continuum (NLMC) method \cite{chung2017non} for problems in heterogeneous fractured media. We remark that since the local solutions are computed in an oversampled domain, the mass transfers between fractures and matrix become non-local, and the resulting upscaled model contains more effective properties of the flow problem.

In this paper, we construct the multicontinuum upscaled models based on NLMC. We construct multiscale basis functions in each local domain for both fractures and matrix. Upscaled model have only one additional coarse degree of freedom (DOF) for each fracture network. Numerical results show that the coupled NLMC method for the fractured porous media provide accurate and efficient upscaled model on the coarse grid.
The implementation is based on the open-source library FEniCS, where we use geometry objects and the FEniCS interface to the linear solvers \cite{logg2009efficient, logg2012automated}.

This paper is organized as follows.
In Section \ref{sec:model}, we consider mathematical model. Next in Section \ref{sec:fine}, we present finite volume fine grid approximation for the EFM and DFM approaches.
In Section \ref{sec:coarse}, we propose an upscaled coarse-grid model for flow in fractured porous media. After that, we present some numerical examples for various model problems, and we show that proposed method can achieve a good accuracy with a very few degrees of freedom and discuss the details in Sections \ref{sec:num1}-\ref{sec:num3}. A conclusion is drawn in Section \ref{sec:conclusion}.

\section{Problem formulation} \label{sec:model}

In this paper, we consider a mixed dimensional mathematical model for the fractured porous medium. This coupled problem describes fluid flow in a domain $\Omega \in \mathcal{R}^d$ (d = 2,3), where the thickness of the fractures are small and can be represented by a reduced dimensional object $\gamma \in \mathcal{R}^{d-1}$. The resulting model reads
\begin{equation}
\label{mm1}
\begin{split}
& \frac{ \partial (\rho \phi_m )}{\partial t} 
+ \nabla \cdot ( \rho q_m) + \rho  r_{mf} = \rho f_m, \quad  x \in \Omega \\
& q_m = -\frac{k_m}{\mu}  \nabla p_m, \quad x \in \Omega, \\
& \frac{ \partial (\rho b)}{\partial t} 
+ \nabla \cdot ( \rho q_f) + \rho  r_{fm} = \rho f_f, \quad 
 x \in \gamma \\
& q_f = - b \frac{k_f}{\mu}  \nabla p_f, \quad  x \in \gamma , 
\end{split}
\end{equation}
where $q_m$ is the velocity in the porous matrix $\Omega$,  $q_f$ is the velocity in the fractures $\gamma$,  $\mu$  is the fluid viscosity, $c_{\alpha}$, $k_{\alpha}$  are the compressibility and permeability ($k_f = k_f(b)$), $f_{\alpha}$ is the source term with $\alpha = f, m$. 

For the coupling term between the fractures and matrix, we have
\[
r_{mf} = \eta_m \sigma (p_m - p_f), \quad
r_{fm} = \eta_f \sigma (p_f - p_m), 
\]
where $\sigma = k^*/b$, $b$ is the fracture thickness and $k^*$ is the harmonic average between $k_m$ and $k_f$. Coefficients $\eta_m$ and $\eta_f$ depends on mesh parameters.
This term expresses the conservation of mass between the two continua.

Let  $\rho = const$, then we have following mixed dimensional coupled system of equation
\begin{equation}
\label{mm2}
\begin{split}
& a_m \frac{ \partial p_m }{\partial t} 
- \nabla \cdot (b_m \nabla p_m) +  \eta_m \sigma (p_m - p_f) =   f_m, 
\quad  x \in \Omega_m, \\
& a_f \frac{ \partial p_f }{\partial t} 
- \nabla \cdot (b_f \nabla p_f) -  \eta_f \sigma (p_m - p_f) =   f_f.
\quad  x \in \gamma,
\end{split}
\end{equation}
where $p_m$ is the pressure in the porous matrix $\Omega$, $p_f$ is the pressure in the fractures $\gamma$, $a_m = c_m$, $a_f = c_f$, $b_m = k_m/\mu$, $b_f = b k_f/\mu$ are constants for simplicity.

\section{Fine-grid finite volume approximation}\label{sec:fine}

Next, we consider fine-grid approximation of the mixed dimensional coupled problem on unstructured grids using cell centered finite volume method.
\begin{figure}[h!]
\centering
\includegraphics[width=0.49 \textwidth]{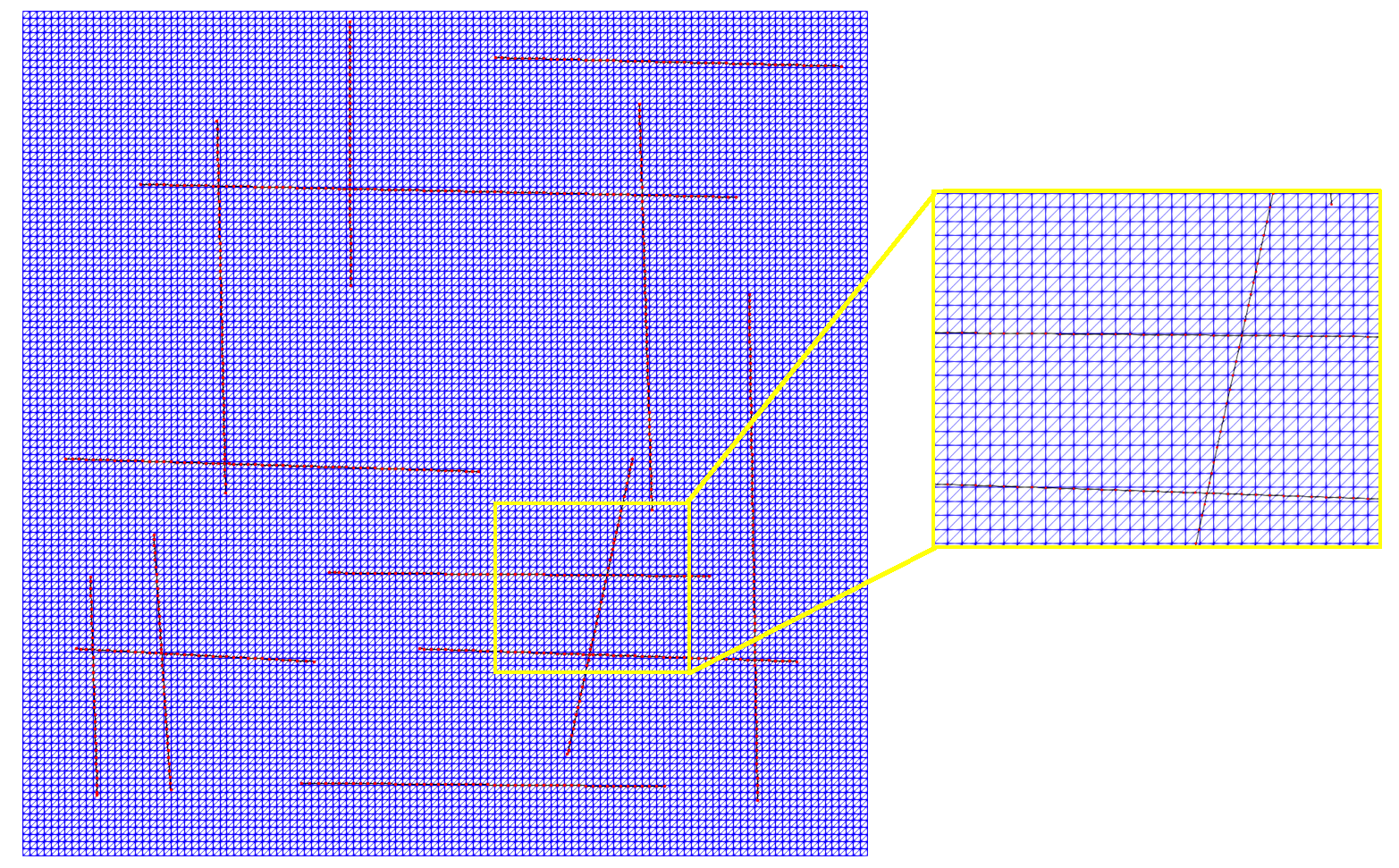}
\includegraphics[width=0.49 \textwidth]{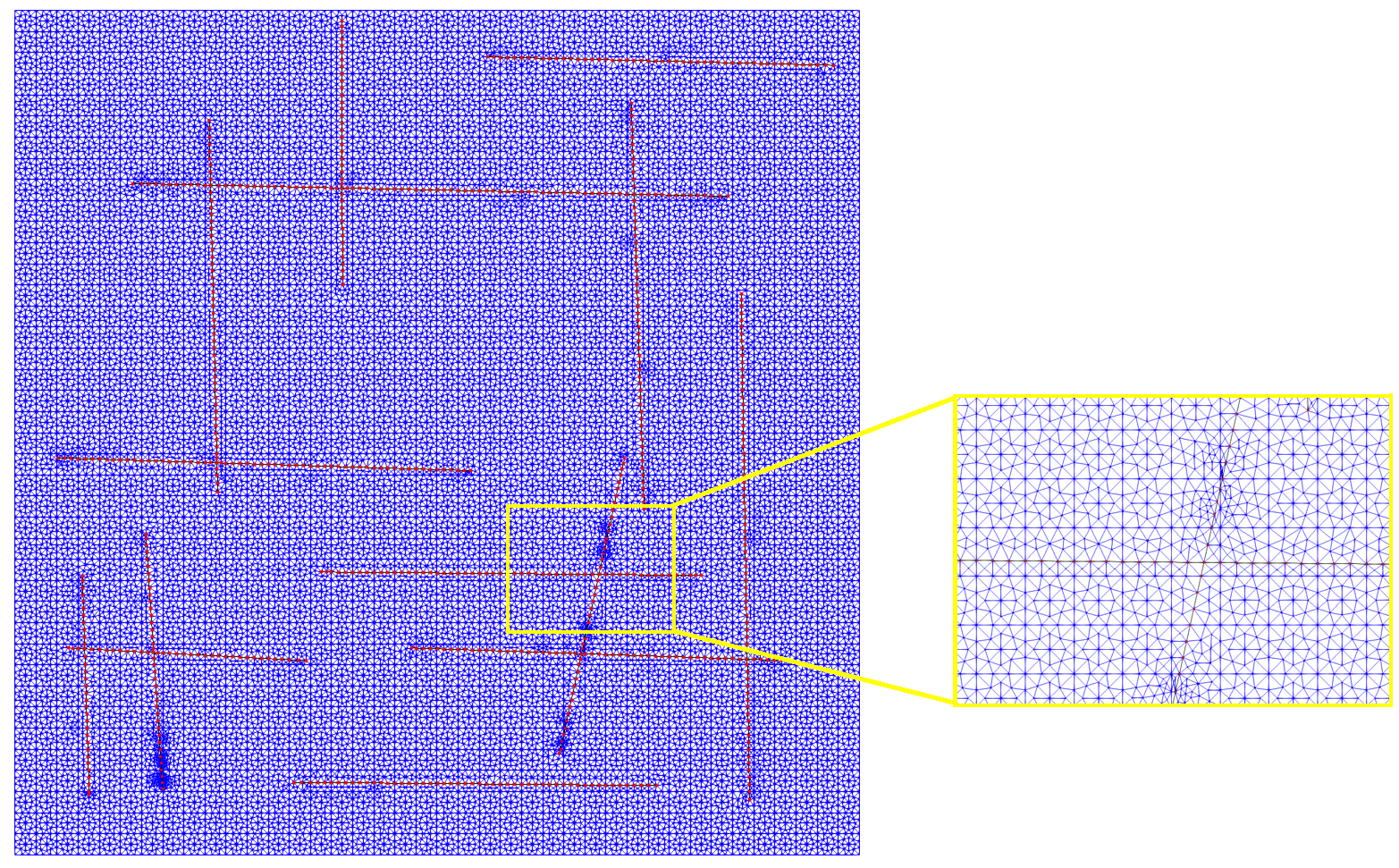}
\caption{Computational grids for EFM (left) and DFM (right).}
\label{fig:sch}
\end{figure}

\textbf{Discrete fracture model approximation. } 
Let $\mathcal{T}_h = \cup_i \varsigma_i$ be the fine mesh of the domain $\Omega$ containing triangular or tetrahedral elements that are conforming with fractures, and let $\mathcal{E}_h$ be all fine-scale facets of the mesh $\mathcal{T}_h$. Denote by $\mathcal{E}_{\gamma} = \cup_l \iota_l$ the fracture facets, where $\mathcal{E}_{\gamma} \subset \mathcal{E}_h$(see right of Figure \ref{fig:sch}).
For discrete fracture model, we have the following discrete problem using two-point flux approximation  
\begin{equation}
\label{dfm1}
\begin{split}
& a_m \frac{ p_{m, i} - \check{p}_{m, i} }{\tau} |\varsigma_i| 
 + \sum_{E_{ij} \subset \partial K_i / \mathcal{E}_{\gamma} }  T_{ij}  (p_{m, i} - p_{m, j})
 +  \sigma_{il} (p_{m, i} - p_{f, l} )  
 =  f_m   |\varsigma_i|, \quad \forall i = 1, N^m_f \\
& a_f \frac{ p_{f, l} - \check{p}_{f, l}}{\tau}  |\iota_l| 
+ \sum_{n} W_{ln} (p_{f, l} - p_{f, n})
- \sigma_{il} (p_{m, i} - p_{f, l} ) 
 =  f_f  |\iota_l|, \quad \forall l = 1, N^f_f
\end{split}
\end{equation}
where 
$T_{ij} = b_m |E_{ij}|/d_{ij}$ ($|E_{ij}|$ is the length of facet between cells $\varsigma_i$ and $\varsigma_j$, $d_{ij}$ is the distance between midpoint of cells $\varsigma_i$ and $\varsigma_j$),  
$W_{ln} = b_f/d_{ln}$ ($d_{ln}$ is the distance between points $l$ and $n$), 
$N^m_f$ is the number of cells in $\mathcal{T}_h$, 
$N^f_f$ is the number of cell related to fracture mesh $\mathcal{E}_{\gamma}$, 
$\sigma_{il} = \sigma$ if $\mathcal{E}_{\gamma} \cap \partial \varsigma_i = \iota_l$ and zero else. 
Here, we take $\eta_m = 1 / |\varsigma_i|$, $\eta_f = 1 / |\iota_l|$ and use implicit scheme for time discretization and $(\check{p}_m, \check{p}_f)$ is the solution from previous time step and $\tau$ is the given time step. 

Note that, the discrete fracture approximation can be used for fluid flow simulation in a fractured porous medium with both high and low permeable fractures.

\textbf{Embedded fracture model approximation. } 
Let $\mathcal{T}_h=  \cup_i \varsigma_i$ be the structured fine grid with triangular or tetrahedral cells of the domain $\Omega$, we note that in this approach the mesh does not need to be conforming with the fracture lines. The additional fracture mesh denoted by $\mathcal{E}_{\gamma} = \cup_l \iota_l$ is only performed on the fractures (see left of Figure \ref{fig:sch}).
For embedded fracture model, we have 
\begin{equation}
\label{efm1}
\begin{split}
& a_m \frac{ p_{m, i} - \check{p}_{m, i} }{\tau} |\varsigma_i| 
 + \sum_{E_{ij} \in \partial K_i}  T_{ij}  (p_{m, i} - p_{m, j})
 +  \sigma_{il} (p_{m, i} - p_{f, l} )  
 =  f_m   |\varsigma_i|, \quad \forall i = 1, N^m_f \\
& a_f \frac{ p_{f, l} - \check{p}_{f, l}}{\tau}  |\iota_l| 
+ \sum_n W_{ln} (p_{f, l} - p_{f, n})
- \sigma_{il} (p_{m, i} - p_{f, l} ) 
 =  f_f  |\iota_l|, \quad \forall l = 1, N^f_f
\end{split}
\end{equation}
where $\sigma_{il} = \sigma$ if $\iota_l \subset \varsigma_i$ and zero else. 

\textbf{Matrix form. }
Therefore, we have the following system of equations for $p = (p_m, p_f)^T$ presented in the matrix form
\begin{equation}
\label{mm-matrix}
M \frac{p -  \check{p} }{\tau} + A p = F,
\end{equation}
\[
M = 
\begin{pmatrix}
M_m & 0 \\
0 & M_f \\
\end{pmatrix}, \quad 
A = 
\begin{pmatrix}
A_m + Q & -Q \\
-Q & A_f+Q
\end{pmatrix}, \quad 
F  =
\begin{pmatrix}
F_m \\
F_f 
\end{pmatrix},
\]
and
\[
M_m = \{m^m_{ij}\}, \quad 
m^m_{ij} = 
\left\{\begin{matrix}
 a_m |\varsigma_i| / \tau & i = j, \\ 
0 & i \neq j
\end{matrix}\right. , \quad 
M_f = \{m^f_{ln}\}, \quad 
m^f_{ln} = 
\left\{\begin{matrix}
 a_f |\iota_l| / \tau & l = n, \\ 
0 & l \neq n
\end{matrix}\right. ,
\]\[
Q = \{q_{il}\}, \quad 
q_{il} = 
\left\{\begin{matrix}
\sigma & i = l, \\ 
0 & i \neq l
\end{matrix}\right. ,
\]
where
$A_m = \{T_{ij}\}$, 
$A_f = \{W_{ln}\}$, 
$F_m = \{f^m_i\}$, $f^m_i = f_m |\varsigma_i|$, 
$F_f = \{f^f_l\}$, $f^f_l = f_f |\iota_l|$ and size of fine-grid system is $N_f = N^m_f + N^f_f$.

\section{Coarse-grid upscaled method}\label{sec:coarse}

Next, we describe the construction of the upscaled model on coarse grid using Non-local multi-continuum (NLMC) approach. In this method, the multiscale basis functions are constructed by solving local problems in the oversampled local region. The basis functions satisfy the constraint that it vanishes in all other continuum except for the target continuum which it is formulated for. Construction of basis is similar for both discrete and embedded fine-grid fracture models.
 
\begin{figure}[h!]
\centering
\includegraphics[width=0.8 \textwidth]{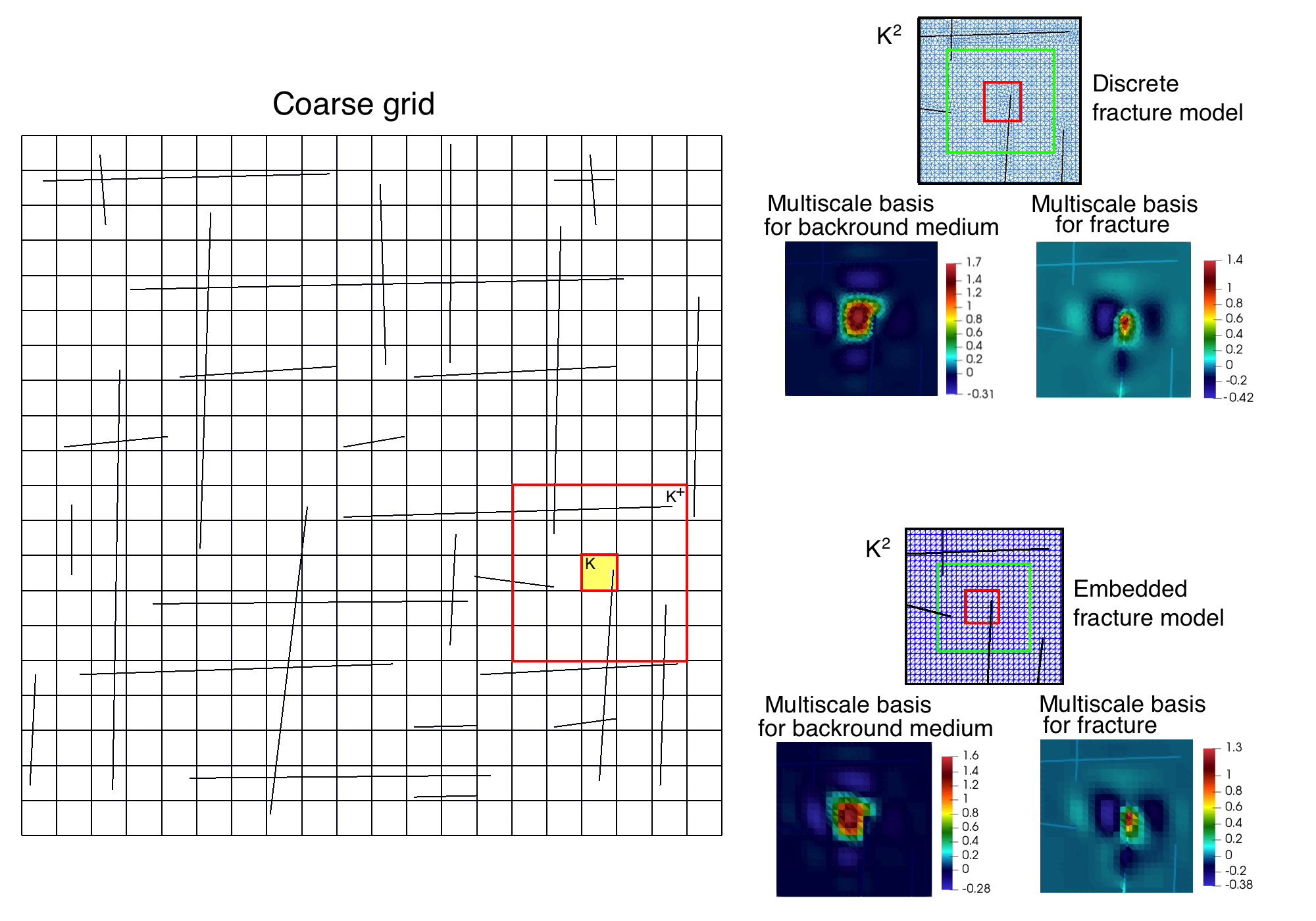}
\caption{Multiscale basis functions on mesh $20 \times 20$ for local domain $K^2$ for DFM and EFM fine-grid approximations }
\label{fig:msbf20}
\end{figure}

In NLMC \cite{chung2017non}, we apply simplified basis for fractured media to form the auxiliary space, which will be used together with an energy minimization principle
for form the required basis functions.  
The resulting multiscale basis functions have spatial decay property in local domains and separate background medium and fractures. 
Finally, the basis functions are used in the construction of the upscaled model. 

Let $K^+_i$ be an oversampled region for the coarse cell $K_i$ (see Figure \ref{fig:msbf20}) obtained by enlarging $K_i$ by several coarse cell layers. 
We will construct a set of basis functions, whose supports are $K_i^+$. 
Each of these basis functions is related to the matrix component in $K_i$ as well as each fracture network within $K_i$. 
For fractures, we denote $\gamma = \cup_{l = 1}^L \gamma^{(l)}$, where $\gamma^{(l)}$ denotes the $l$-th fracture network and $L$ is the total number of fracture networks.
We also write $\gamma^{(l)}_j = K_j \cap \gamma^{(l)}$ is the fracture  inside coarse cell $K_j$ and  $L_j$ is the number of fracture networks in $K_j$. 
For each $K_i$, we will therefore obtain $L_j+1$ basis functions: one for $K_i$ and one for each $\gamma^{(l)}_i$. 
Following the framework of \cite{chung2017non} and \cite{chung2017constraint},
the auxiliary space $V^{aux}(K_i)$ for the coarse cell $K_i$ contains functions that are supported in $K_i$
and are piecewise constant functions such that they are constant on $K_i$ and on each $\gamma^{(l)}_i$. 


We next define the constraints that will be used for multiscale basis construction. 
For each $K_j \subset K_i^+$: \\
(1) background medium ($\psi_0^{i}$) :
\[
\int_{K_j} \psi_0^i \, dx = \delta_{i,j}, \quad 
\int_{\gamma^{(l)}_j} \psi_0^i \, ds = 0, \quad  l=\overline{1, L_j}.
\]
(2) $l$-th fracture network in $K_i$ ($\psi_l^i$):
\[
\int_{K_j} \psi_l^i \, dx = 0, \quad 
\int_{\gamma^{(l)}_j} \psi_l^i \, ds = \delta_{i,j}\delta_{m,l}, \quad  l=\overline{1, L_j}.
\]
We first discuss the constraint for background medium in (1). We note that it is a set of constraints
so that the resulting function has mean value one on the coarse cell $K_i$, and has mean value zero on all other coarse cells within $K_i^+$.
In addition, the resulting function has mean value on all fracture networks within $K_i^+$. 
We next discuss the constraint for the fracture network (2). 
We note that it is a set of constraints
so that the resulting function has mean value zero on all coarse cells within $K_i^+$.
Moreover, the resulting function has mean value one on the target fracture network $\gamma^{(l)}_i$ and has mean value zero
on all fracture networks within $K_i^+$. 
To sum up, the above constraints will give $L_i+1$ functions. 

Together with the above constraints, we will construct the basis functions as follows. 
Following the framework of \cite{chung2017non} and \cite{chung2017constraint}, we will find the multiscale basis functions using the energy minimizing constraint property. As a result, we will solve the following local problems in the oversampled region $K_i^+$ using a fine-grid approximation for the system of flow in fractured porous media presented in the previous Section. In particular, we solve following coupled system in  $K_i^+$:
\begin{equation}
\label{eq:basis}
\begin{pmatrix}
A_m+ Q & -Q & B^T_m & 0 \\
-Q & A_f + Q & 0 &  B^T_f \\
B_m & 0 & 0 & 0 \\
0 & B_f & 0 & 0 \\
\end{pmatrix} 
\begin{pmatrix}
\psi_m \\
\psi_f \\
\mu_m \\
\mu_f \\
\end{pmatrix} = 
\begin{pmatrix}
0 \\
0 \\
F_m \\
F_f \\
\end{pmatrix}
\end{equation}
with  zero Dirichlet boundary conditions on $\partial K^+_i$ for $\psi_m$ and $\psi_f$. 
We remark that $(\psi_m,\psi_f)$ denotes each of the basis functions that satisfy the above constraints. 
Note that we used Lagrange multipliers $\mu_m$ and $\mu_f$ to impose the constraints in the multiscale basis construction. 

To construct multiscale basis function with respect to porous matrix $\psi^0 = (\psi^0_m, \psi^0_f)$, we set $F_m = \delta_{i,j}$ and $F_f = 0$. For multiscale basis function with respect to the $l$-th fracture network, we set  $F_m =  0$ and $F_f = \delta_{i,j}\delta_{m,l}$. In Figure \ref{fig:msbf20}, we depict a multiscale basis functions for oversampled region $K^+_i = K^2_i$ (two oversampling coarse cell layers) in a $20 \times 20$ coarse mesh. Combining these multiscale basis functions, we obtain the following multiscale space
\[
V_{ms} = \text{span} \{ (\psi^{i,l}_m, \psi^{i,l}_f), \, i = \overline{1,N_c}, \, l = \overline{0, L_i} \}
\]
and the projection matrix
\[
R = \begin{pmatrix}
R_{mm} & R_{mf} \\
R_{fm} & R_{ff} \\
\end{pmatrix}, 
\]\[
R_{mm}^T = \left[ 
\psi^{0,0}_m, \psi^{1,0}_m  \ldots \psi^{N_c,0}_m
\right],\quad
 R_{ff}^T = \left[ 
\psi^{0,1}_f \ldots \psi^{0,L_0}_f, 
\psi^{1,1}_f \ldots \psi^{1,L_1}_f,
\ldots, 
\psi^{N_c,1}_f  \ldots \psi^{N_c,L_{N_c}}_f
 \right],
\]\[
R_{mf}^T = \left[ 
\psi^{0,0}_f, \psi^{1,0}_f \ldots \psi^{N_c,0}_f
 \right], \quad
R_{fm}^T = \left[ 
\psi^{0,1}_m \ldots \psi^{0,L_0}_m, 
\psi^{1,1}_m \ldots \psi^{1,L_1}_m,
\ldots, 
\psi^{N_c,1}_m  \ldots \psi^{N_c,L_{N_c}}_m
 \right], 
\]

Therefore, the resulting upscaled coarse grid model reads
\begin{equation}
\label{t-nlmc2}
\bar{M} \frac{\bar{p}^{n+1} - \bar{p}^n}{\tau} + \bar{A} \bar{p}^{n+1} = \bar{F},
\end{equation}
where 
$\bar{A} = R A R^T$, $\bar{F} = R F$ and
$\bar{p} = (\bar{p}_m, \bar{p}_f)$.
We remark that $\bar{p}_m$ and $\bar{p}_f$ are the average cell solution on coarse grid element for background matrix and for fracture media.
That is, each component of $\bar{p}_m$ corresponds to the mean value of the solution on each coarse cell. 
Moreover, each component of $\bar{p}_f$ corresponds to the mean value of the solution on each fracture network with a coarse cell. 

As an approximation, we can use diagonal mass matrix directly calculated on the coarse grid 
\[
\bar{M} = 
\begin{pmatrix}
\bar{M}_m & 0 \\
0 & \bar{M}_f \\
\end{pmatrix}, \quad
\bar{Q} = 
\begin{pmatrix}
\tilde{Q} &  - \tilde{Q} \\
-\tilde{Q} & \tilde{Q} \\
\end{pmatrix}, \quad
\bar{F} = 
\begin{pmatrix}
\bar{F}_m \\
\bar{F}_f \\
\end{pmatrix}, 
\] 
where 
$\bar{M}_m = \text{diag}\{ a_m |K_i| \}$, 
$\bar{M}_f = \text{diag}\{ a_f |\gamma_i| \}$,  
$\tilde{Q} = \text{diag}\{ \sigma |\gamma_i| \}$ 
and for the right-hand side vector 
$\bar{F}_m = \{ f_m |K_i| \}$, 
$\bar{F}_f = \{ f_f |\gamma_i| \}$. 
We remark that the matrix $A$ is non-local and provide good approximation due to the basis construction.

\section{Numerical results for high and low permeable fractures with DFM}
\label{sec:num1}

\begin{figure}[h!]
\centering
\includegraphics[width=0.6 \textwidth]{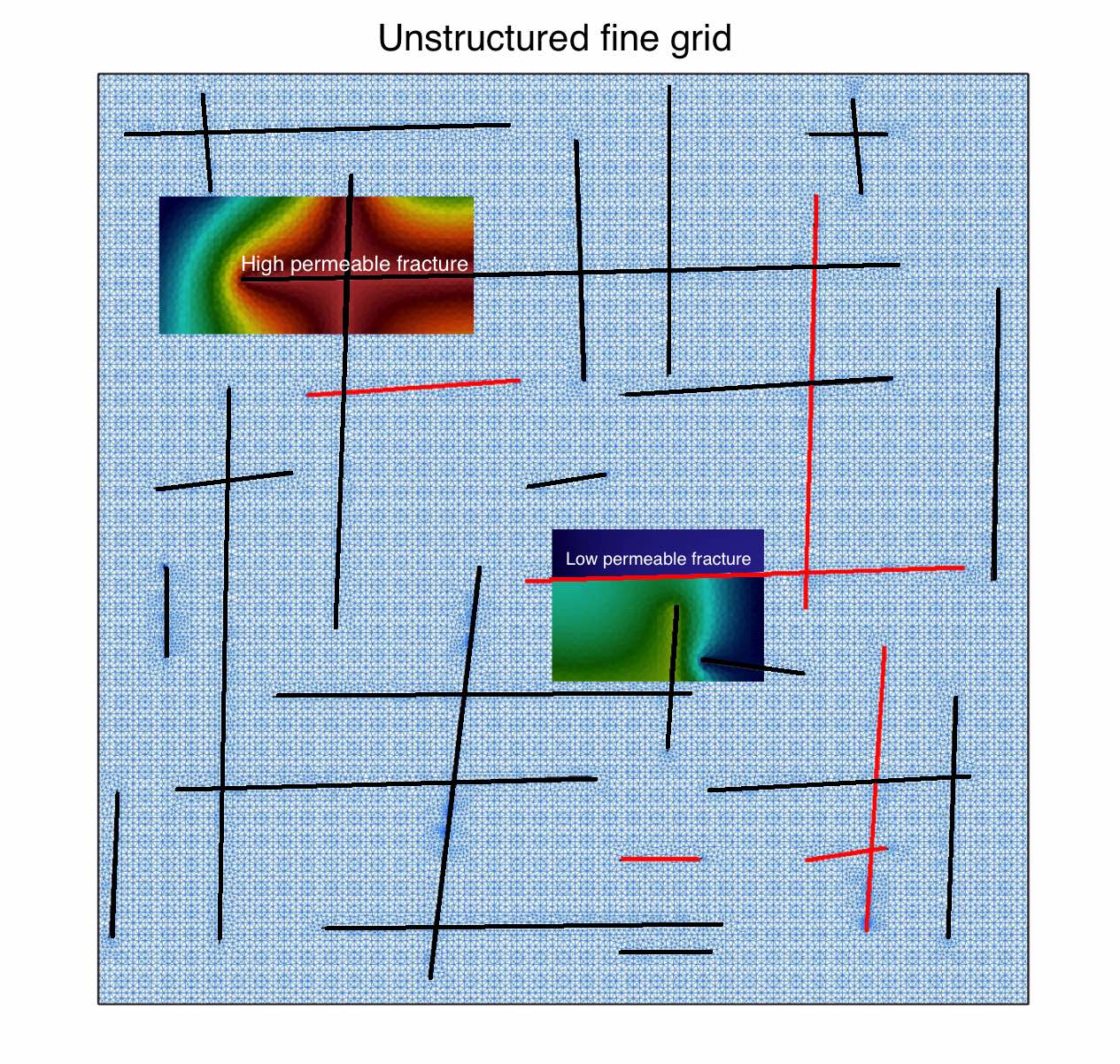}
\caption{Discrete fracture model with low and high permeable fractures. Blue: fine grid. Black: high permeable fractures. Red: low permeable fractures.}
\label{fig:sch-hl}
\end{figure}

In this section, we consider low and high permeable fractures (see Figure \ref{fig:sch-hl} for illustration).  
We construct an accurate approximation of the pressure equation using NLMC approach. The ideas to construct the basis are similar for low and high permeable fracture cases, where only the underlying fine grid models are different. 
In Figure \ref{fig:msbf20}, we depict multiscale basis functions in an oversampled local domain $K^2$ on a $20 \times 20$ coarse mesh.

\begin{figure}[h!]
\centering
Multiscale basis functions for high permeable fracture\\
\includegraphics[width=0.32 \textwidth]{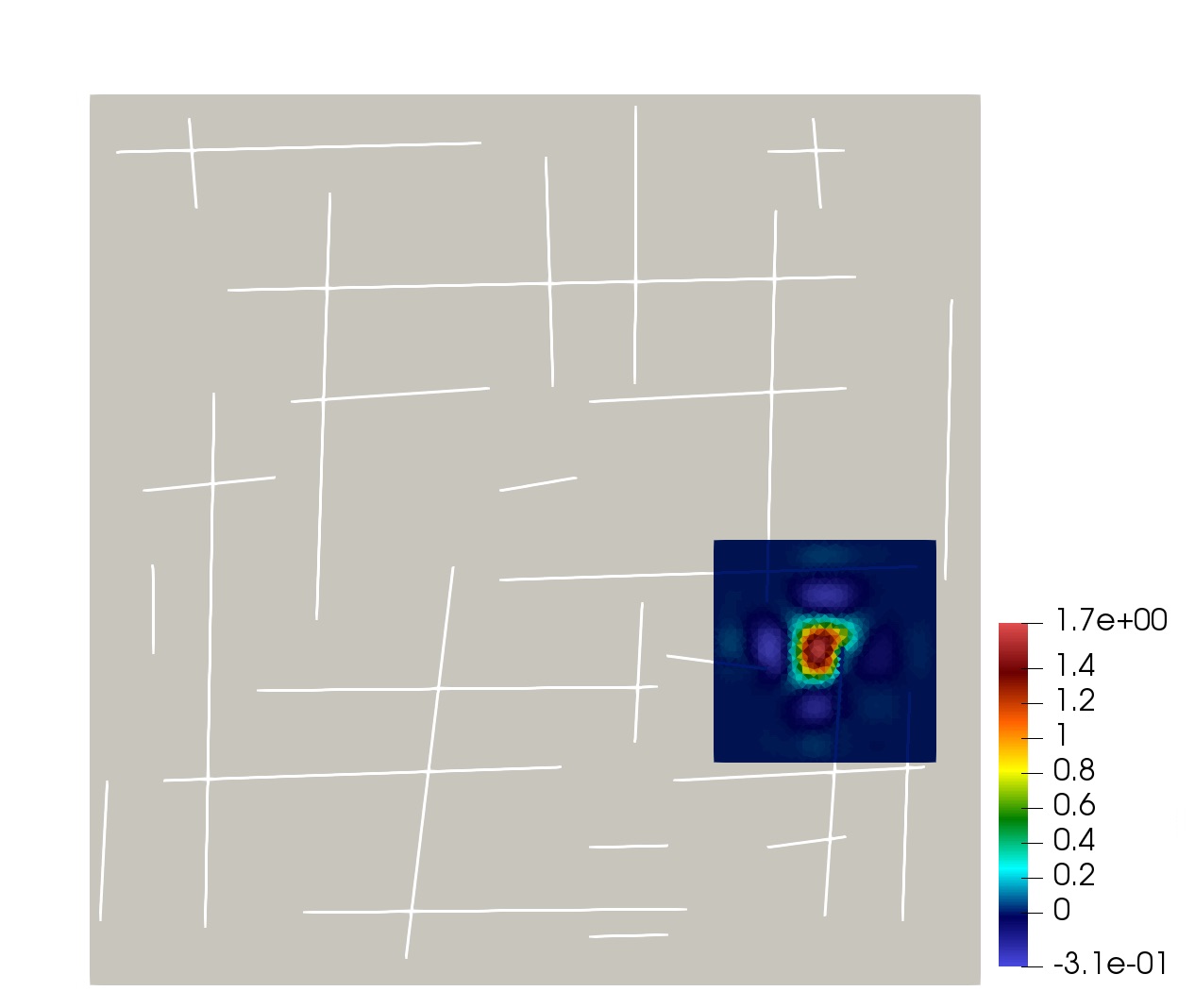}
\includegraphics[width=0.32 \textwidth]{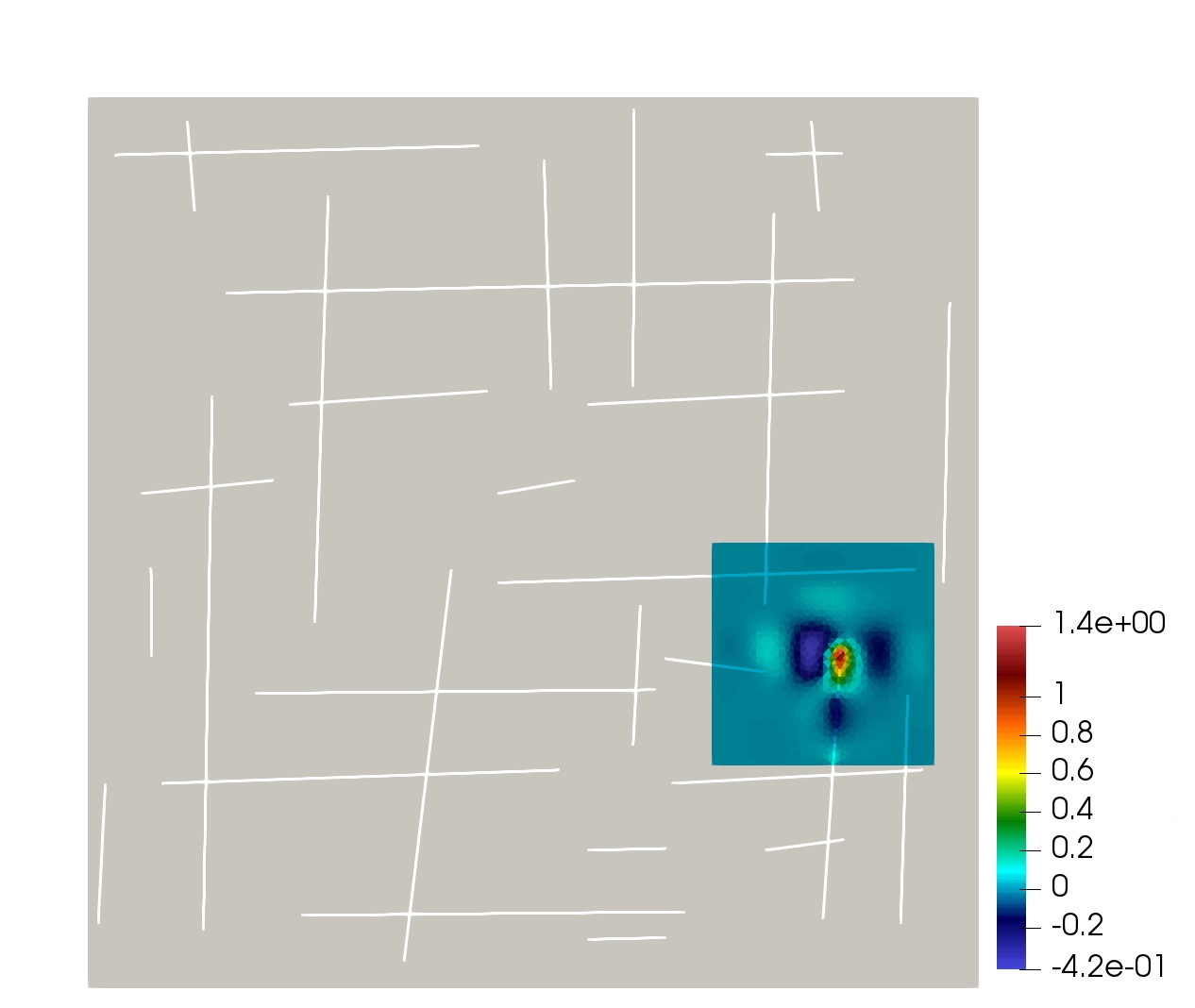}\\
Multiscale basis functions for low permeable fracture\\
\includegraphics[width=0.32 \textwidth]{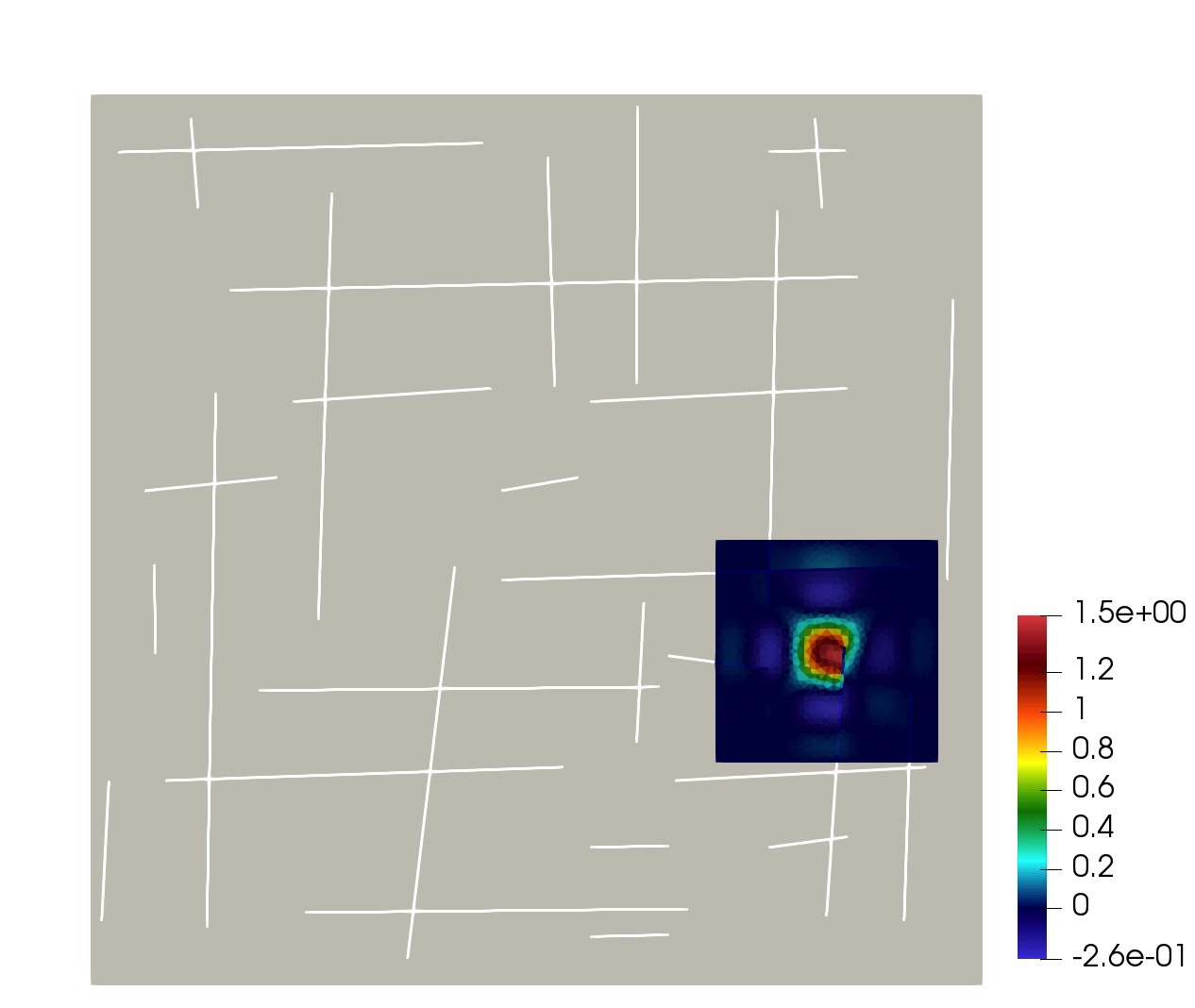}
\includegraphics[width=0.32 \textwidth]{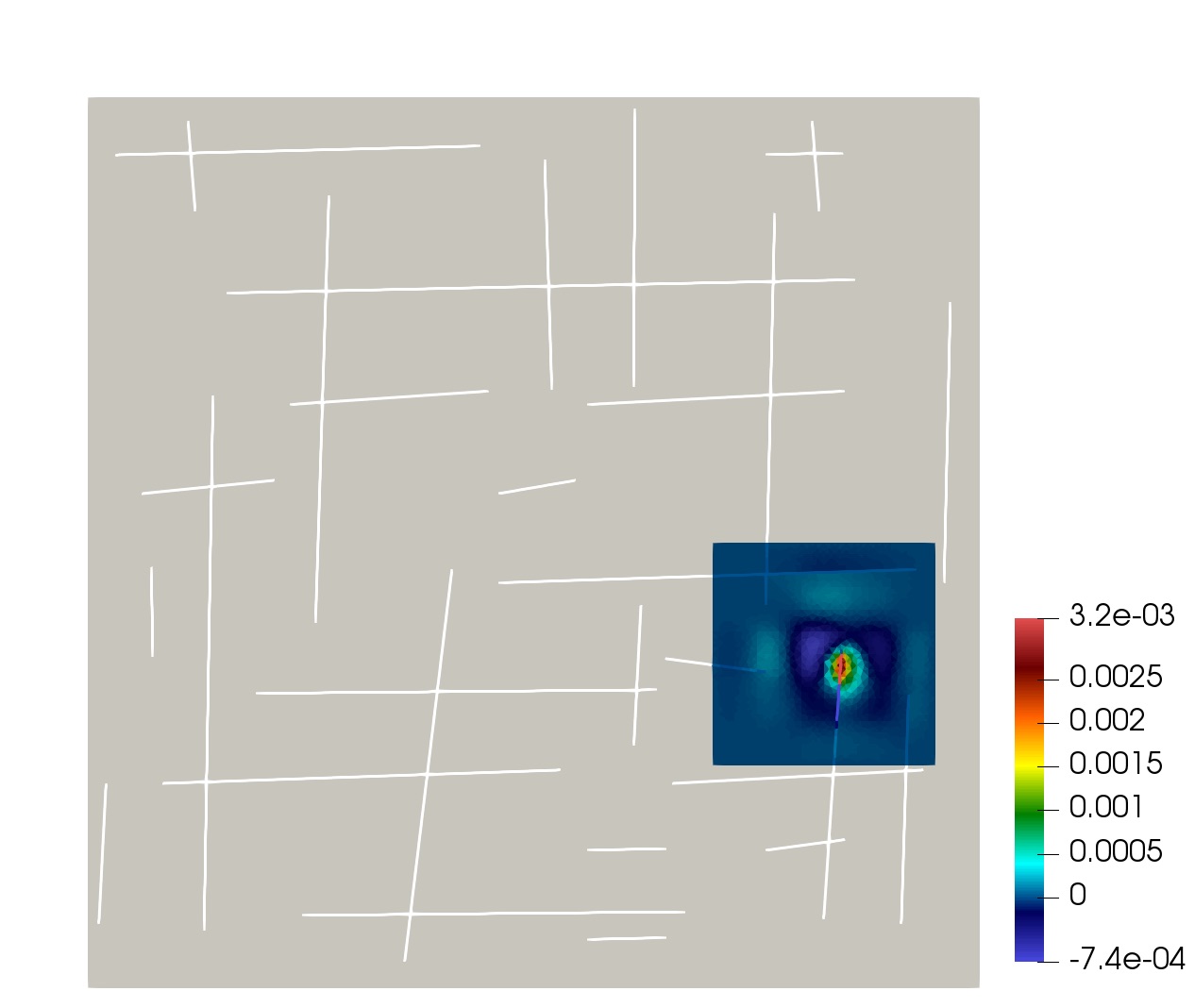}
\caption{Multiscale basis functions on a $20 \times 20$ mesh for local domain $K^2$. First row: multiscale basis functions for matrix and high permeable fracture.
Second row: multiscale basis functions for matrix and low permeable fracture. 
Left column:  multiscale basis functions for matrix.
Right column:  multiscale basis functions for fracture.
}
\label{fig:msbf202}
\end{figure}

\begin{figure}[h!]
\centering
\includegraphics[width=0.32 \textwidth]{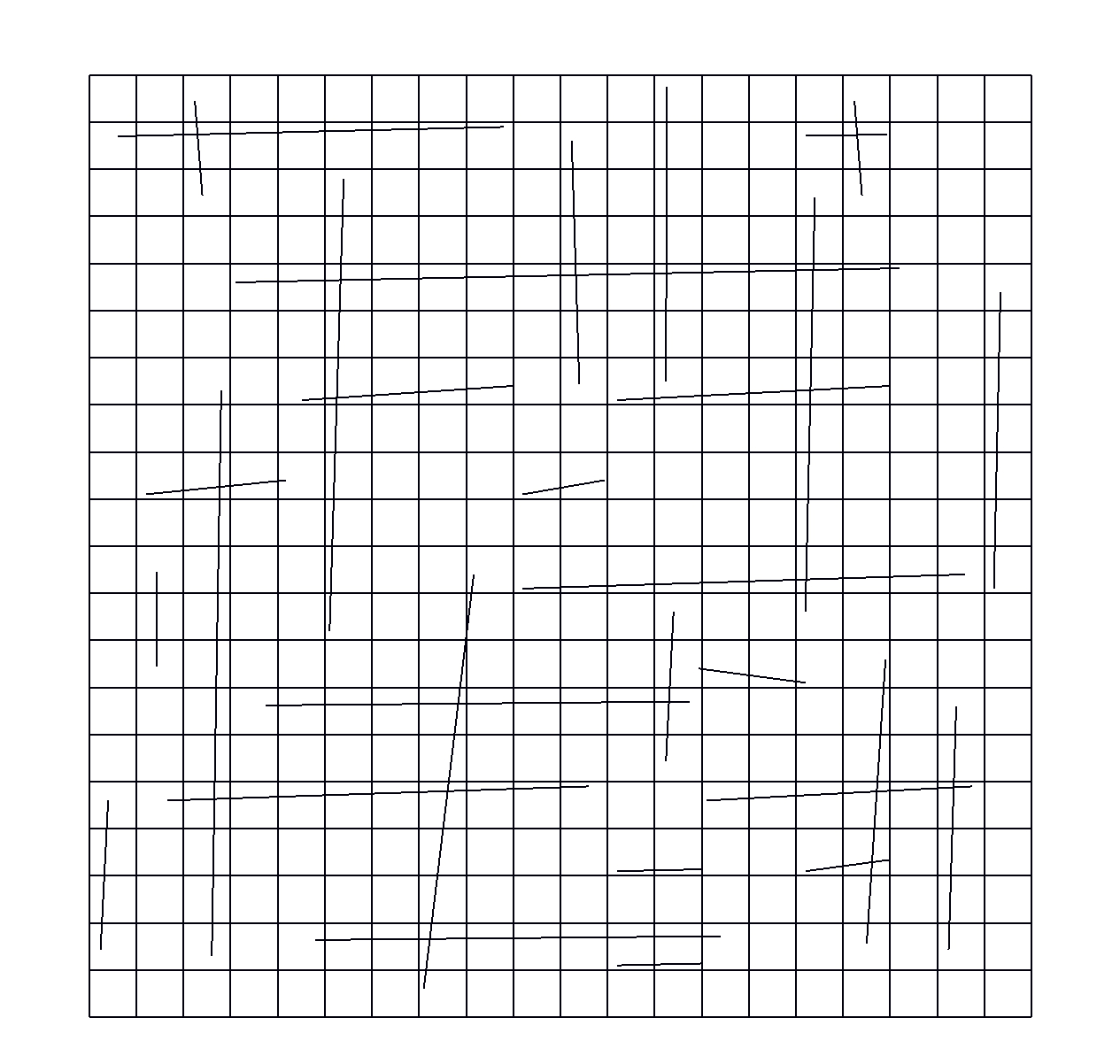}
\includegraphics[width=0.32 \textwidth]{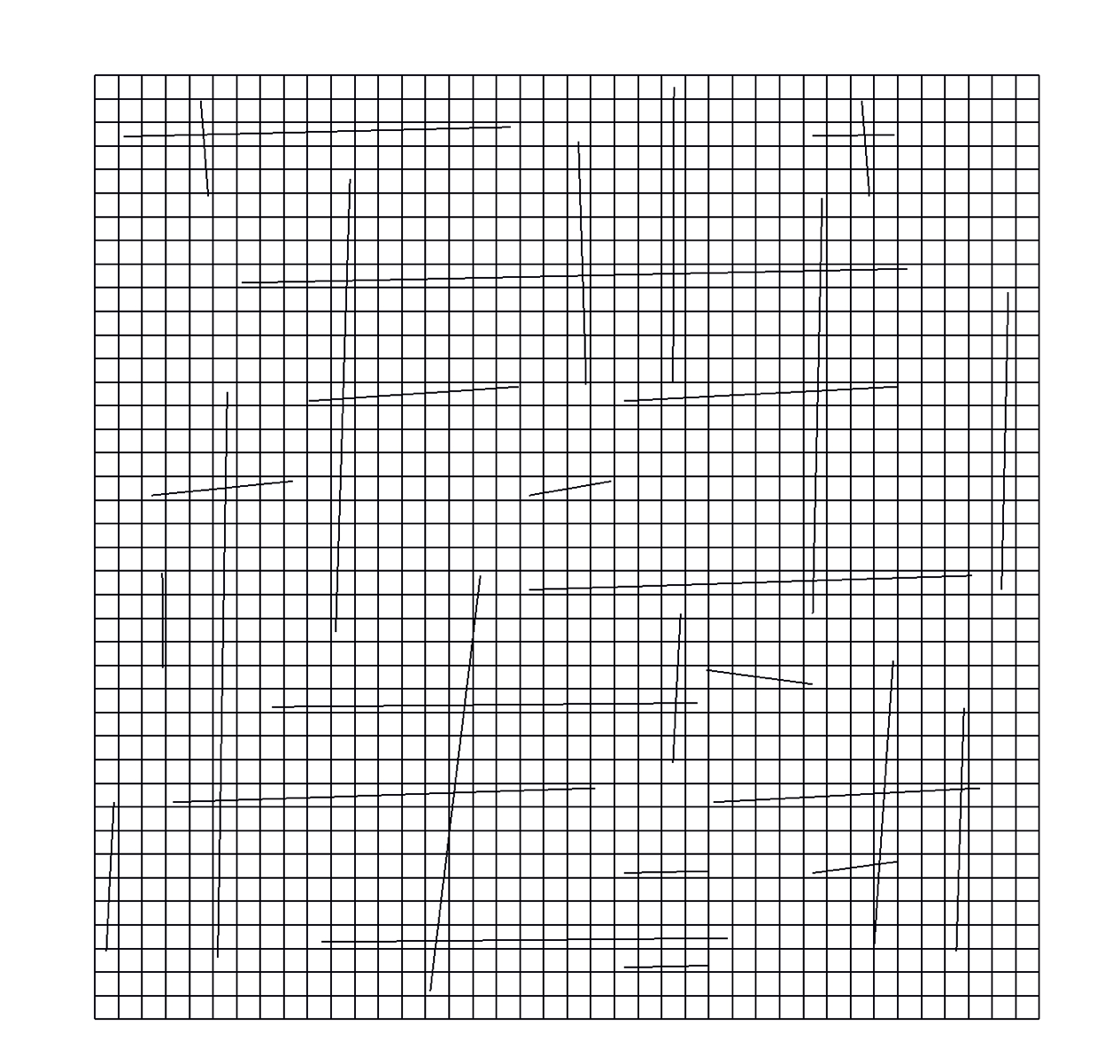}
\includegraphics[width=0.32 \textwidth]{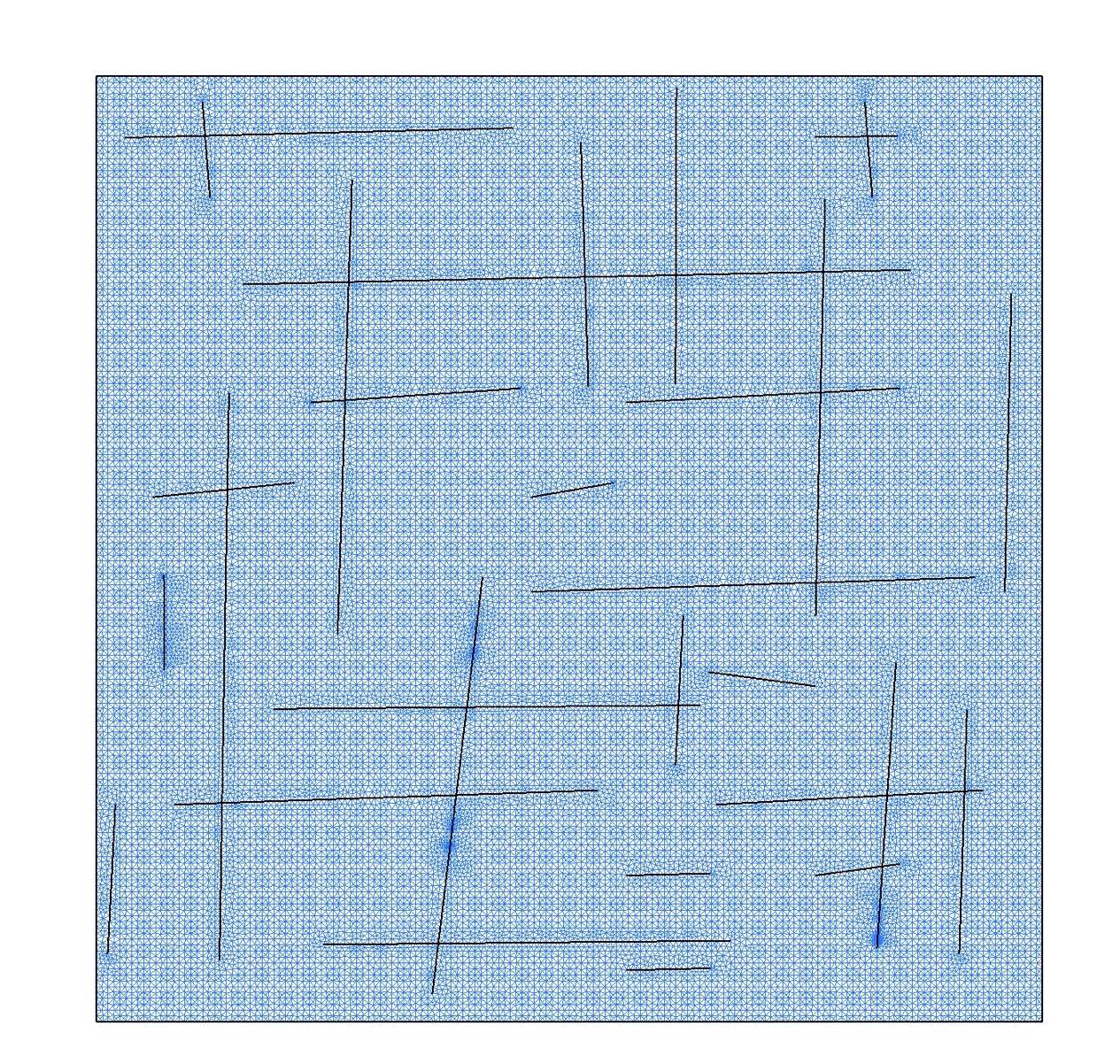}
\caption{Computational grids.  Coarse grids with $400$ and $1600$ cells.
Fine grid with 47520 elements (matrix) and 1042 elements (fractures). }
\label{fig:mesh-hl}
\end{figure}

We consider the computational domain
$\Omega = [0, 1] \times [0, 1]$ with 30 fractures. 
In Figure \ref{fig:mesh-hl}, we show the coarse and fine grids. 
For fine-grid models, we use DFM, thus the fractures are resolved by the fine grid.
The coarse grids are uniformly partitioned into $20 \times 20$ and $40 \times 40$ coarse blocks.  

\begin{figure}[h!]
\centering
\includegraphics[width=1 \textwidth]{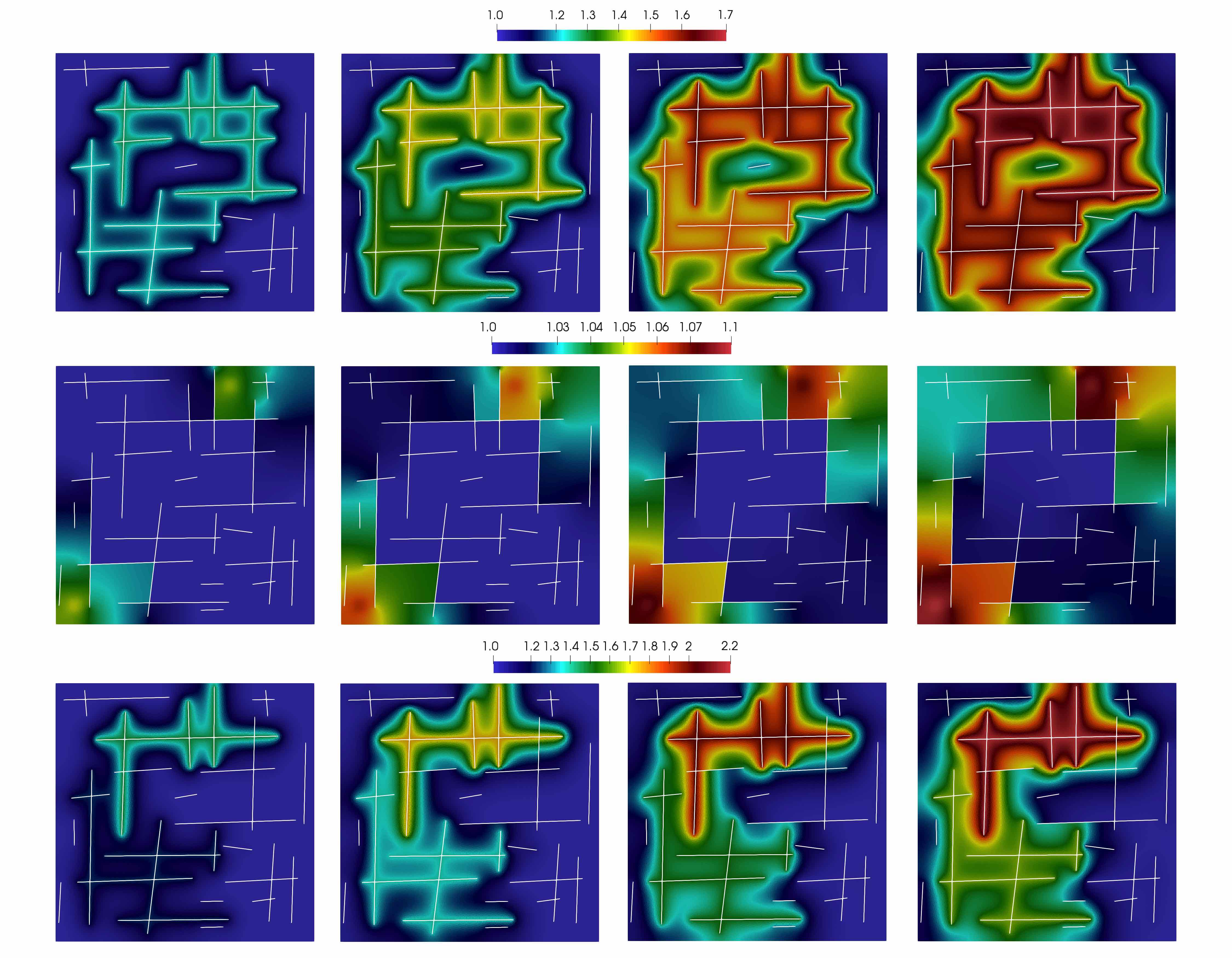}
\caption{Fine-scale solutions for Tests 1,2 and 3 (from top to bottom) for different time steps $t_5 = 0.025$ (first column), $t_{10} = 0.05$ (second column), $t_{15} = 0.075$ (third column) and $t_{20} = 0.1$ (fourth column).  $DOF_f = 48562$. }
\label{fig:uf}
\end{figure}

\begin{figure}[h!]
\centering
\includegraphics[width=1 \textwidth]{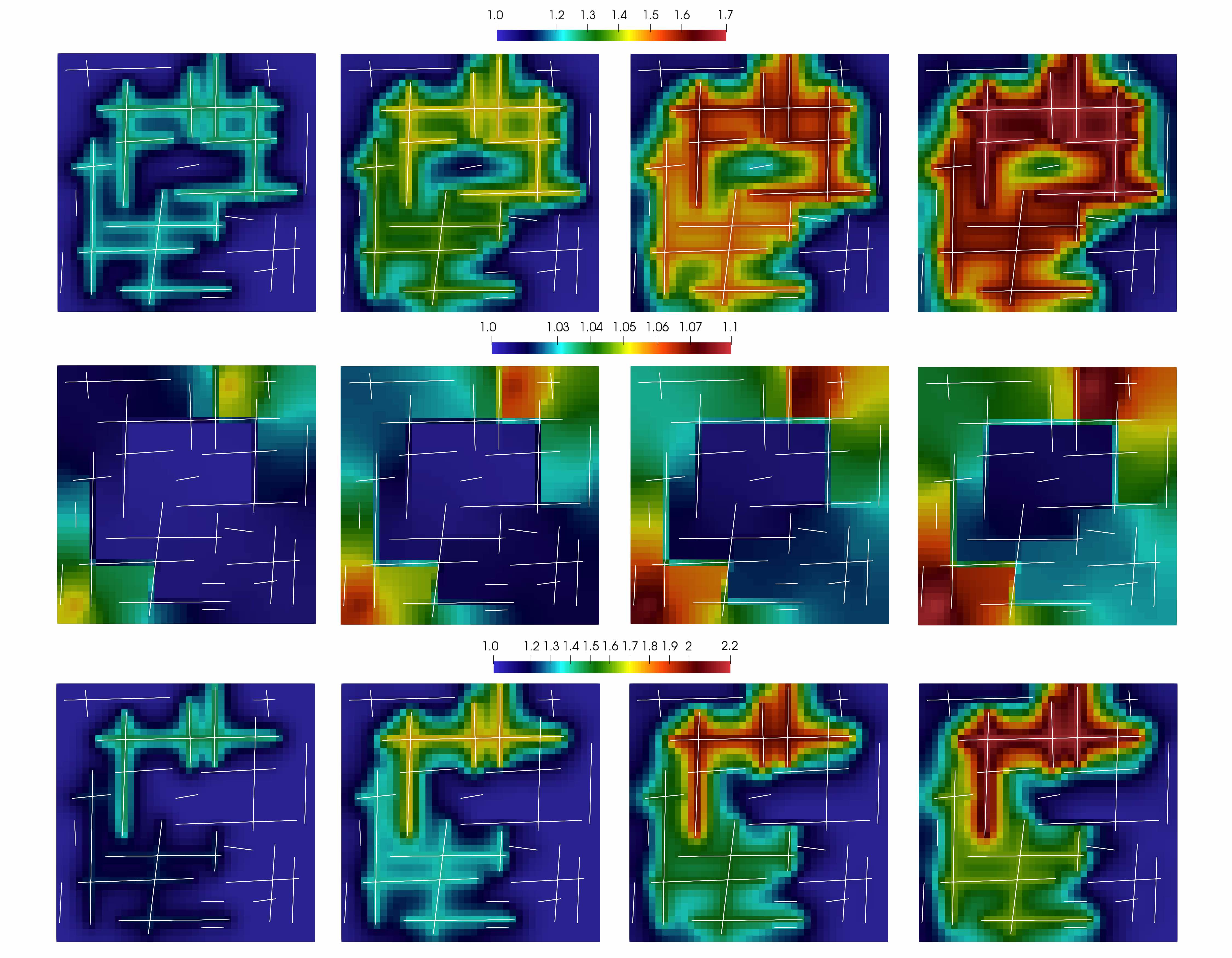}
\caption{Upscaled coarse grid solutions on mesh $40 \times 40$ with $K^3$ for different time steps $t_5 = 0.025$ (first column), $t_{10} = 0.05$ (second column), $t_{15} = 0.075$ (third column) and $t_{20} = 0.1$ (fourth column). Tests 1,2 and 3 (from top to bottom). $DOF_c = 1965$. }
\label{fig:u40abc}
\end{figure}

\begin{table}[h!]
\centering
\begin{tabular}{ |c | c | c | c | c | }
\hline
$K^s$ & 
$t_5$ & 
$t_{10}$ & 
$t_{15}$ & 
$t_{20}$ \\ \hline
\multicolumn{5}{|c|}{Test 1. $20 \times 20$} \\
\hline
$s = 1$	& 0.569	&	1.147	&	1.740	&	2.307 \\ \hline
$s = 2$	& 0.233	&	0.246	&	0.263	&	0.283 \\ \hline
$s = 3$	& 0.229	&	0.236	&	0.254	&	0.272 \\ \hline
\multicolumn{5}{|c|}{Test 2. $20 \times 20$} \\
\hline
$s = 2$	& 0.862	&	2.780	&	6.183	&	11.990 \\ \hline
$s = 3$	& 0.095	&	0.151	&	0.191	&	0.224 \\ \hline
\multicolumn{5}{|c|}{Test 3. $20 \times 20$} \\
\hline
$s = 1$	& 0.793	&	1.891	&	3.108	&	4.319 \\ \hline
$s = 2$	& 0.323	&	0.395	&	0.432	&	0.451 \\ \hline
$s = 3$	& 0.311	&	0.382	&	0.427	&	0.446 \\ \hline
\end{tabular} \,\,\,\,
\begin{tabular}{ |c | c | c | c | c | }
\hline
$K^s$ & 
$t_5$ & 
$t_{10}$ & 
$t_{15}$ & 
$t_{20}$ \\ \hline
\multicolumn{5}{|c|}{Test 1. $40 \times 40$} \\
\hline
$s = 1$	& 0.600	&	1.164	&	1.695	&	2.175 \\ \hline
$s = 2$	& 0.162	&	0.255	&	0.326	&	0.383 \\ \hline
$s = 3$	& 0.151	&	0.202	&	0.231	&	0.248 \\ \hline
\multicolumn{5}{|c|}{Test 2. $40 \times 40$} \\
\hline
$s = 2$	& 0.429	&	1.118	&	2.068	&	3.271 \\ \hline
$s = 3$	& 0.072	&	0.122	&	0.162	&	0.195 \\ \hline
\multicolumn{5}{|c|}{Test 3. $40 \times 40$} \\
\hline
$s = 1$	& 0.839	&	1.751	&	2.699	&	3.610 \\ \hline
$s = 2$	& 0.215	&	0.348	&	0.442	&	0.523 \\ \hline
$s = 3$	& 0.207	&	0.331	&	0.419	&	0.476 \\ \hline
\end{tabular}
\caption{Relative errors of the mean solution on a coarse mesh $20 \times 20$ (Left) and $40 \times 40$ (Right). Test 1, 2 and 3 }
\label{err-abc}
\end{table}

We consider three test cases:
\begin{itemize}
\item Test 1 (high permeable fractures). $k_m = 10^{-6}$, $k_f = 1.0$ for all  fractures.
\item Test 2  (low permeable fractures). $k_m = 10^{-4}$, $k_f = 10^{-10}$ for all  fractures.
\item Test 3 (hybrid fractures). $k_m = 10^{-6}$, $k_f = 1.0$ for 24 fractures and $k_f = 10^{-12}$ for 6 fractures.
\end{itemize}
The other model parameters are chosen as follows:  $c_m = 10^{-5}$, $c_f = 10^{-6}$  with $\sigma =\frac{2.0}{k_m^{-1} + k_f^{-1}}$. 

We set a source $q = 10^{-3}$ on the fractures in the two coarse cells,
\begin{itemize}
\item Test 1 and 3. cell: $0.1<x<0.15, 0.05<y<0.1$ and cell: $0.6<x<0.65, 0.9<y<0.95$.
\item Test 2. cell: $0.05<x<0.1, 0.05<y<0.1$ and cell: $0.65<x<0.7, 0.9<y<0.95$.
\end{itemize}
As for initial pressure, we set $p_0 = 1$. 
Our total simulation time is $t_{max} = 0.1$, and we take 20 time steps for upscaled and fine-scale solvers. 

To compare the results, we investigate the relative $L^2$ error between coarse cell average of the fine-scale solution $\bar{p}_f$ and upscaled coarse grid solutions $\bar{p}$
\begin{equation}
e_{L_2} = ||\bar{p}_f - \bar{p} ||_{L^2}, \quad
|| \bar{p}_f - \bar{p} ||^2_{L^2} =  
\frac{ \sum_K (\bar{p}^K_f - \bar{p}^K)^2 \, dx}{\sum_K (\bar{p}^K_f)^2 \, dx}, \quad 
\bar{p}^K_f = \frac{1}{|K|} \int_K p_f \, dx.
\end{equation}

The fine-scale system has dimension $DOF_f = 47520 + 1042$. The upscaled model has dimension $DOF_c = 593$ for coarse mesh with 400 cells ($20 \times 20$, and $DOF_c = 1965$ for coarse mesh with 1600 cells ($40 \times 40$).
In Figure \ref{fig:uf}, we present the fine scale solution for all test cases at different time steps $t_5 = 0.025$, $t_{10} = 0.05$, $t_{15} = 0.075$ and $t_{20} = 0.1$ . The first row present solutions for Test 1, where we have highly conductive fractures. The second row show the solutions for the low permeable fractures. Finally, in the third row we depict solutions where the fractures have both high and low permeability.

In Figure \ref{fig:u40abc}, we present the upscaled solutions for coarse grid  $40 \times 40$ for Test 1, 2 and 3. For basis calculations, we use oversampled domain $K^+$ with 3 coarse cells layers oversampling. 
We observe  good accuracy of the proposed method with less than one percent of error for all test cases. 

In Table \ref{err-abc}, we present relative errors for two coarse grids and for different numbers of oversampling layers $K^s$ with $s = 1,2$ and $3$. 
We notice a huge reduction of the system dimension and very small errors for unsteady mixed dimensional coupled system.

\section{Numerical results with EFM} 
\label{sec:numerical}

\begin{figure}[h!]
\centering
\includegraphics[width=0.49 \textwidth]{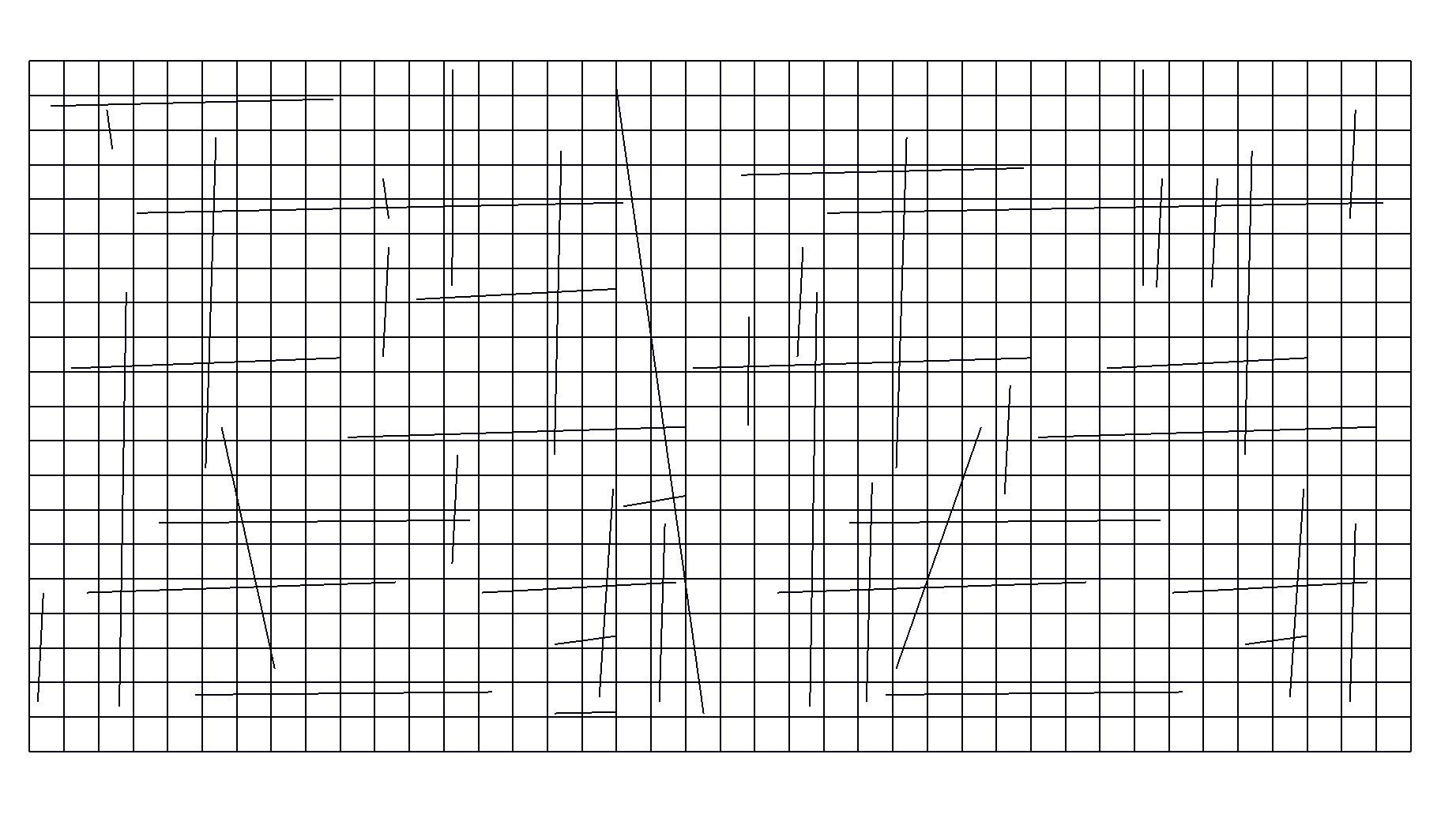}
\includegraphics[width=0.49 \textwidth]{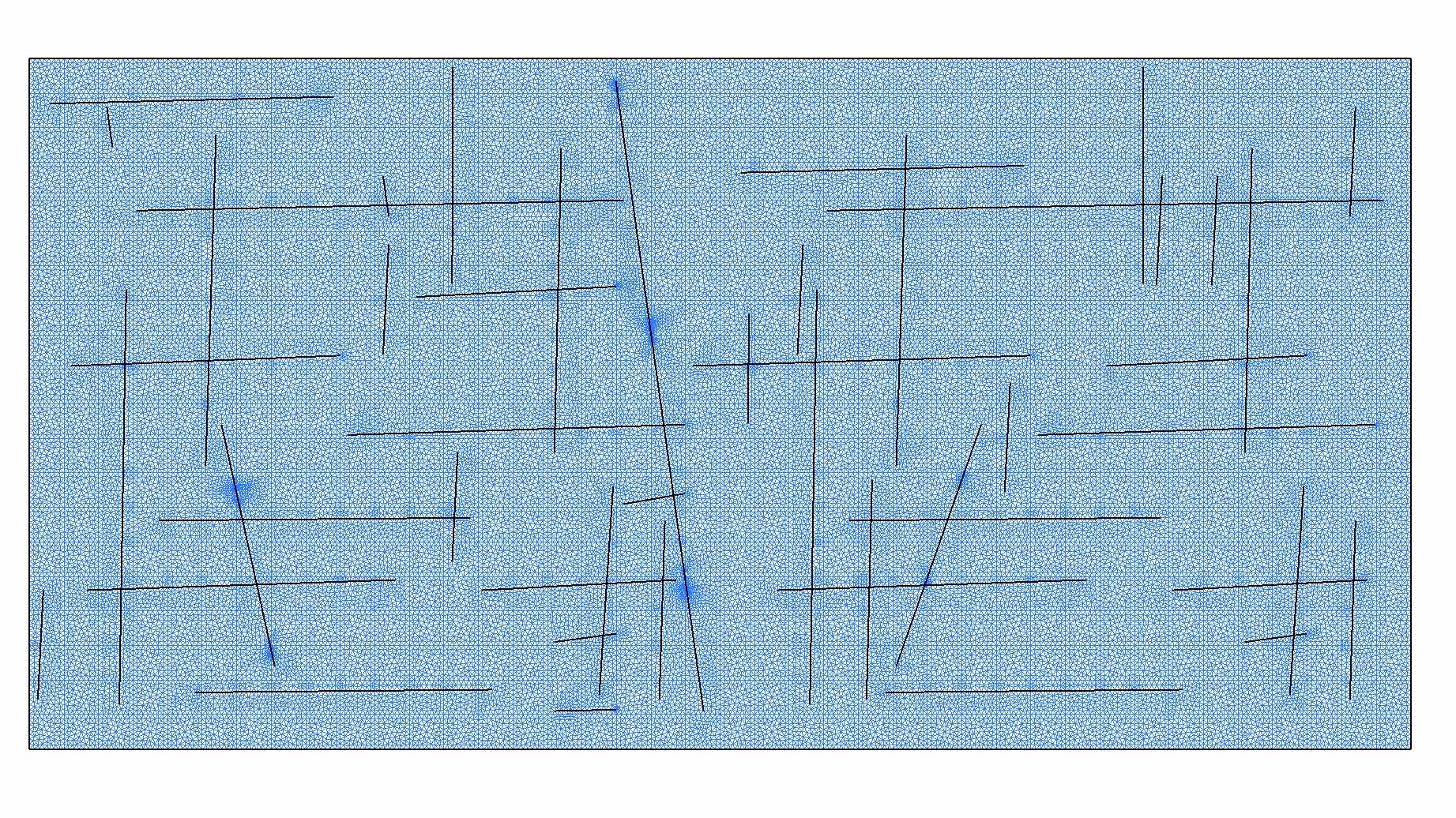}
\caption{Computational grids.  Coarse grid with 800 cells.
Fine grid with 98398 cells (matrix) and 2170 cells (fractures) (DFM). \textit{Geometry 2}
}
\label{fig:mesh2}
\end{figure}

In this section, we present numerical results for upscaled model for embedded and discrete fine-grid fracture models. We consider highly permeable fractures for two types of geometries. 
As for Geometry 1, we consider 30 fracture lines in the domain (Test 1 from the previous section) with injection and production wells (Figure \ref{fig:mesh-hl}). 
Geometry 2 is the computational domain $\Omega = [0, 2] \times [0, 1]$ with 50 fractures (Figure \ref{fig:mesh2}). 

We set a source term on the fractures inside following cells:
\begin{itemize}
\item Geometry 1. Cell $0.1<x<0.15$, $0.05<y<0.1$ (injection) and cell $0.6<x<0.65$, $0.9<y<0.95$ (production) with $q = \pm 10^{-3}$.
\item Geometry 2. Cells $0.1<x<0.15$, $0.05<y<0.1$ and $1.6<x<1.65$, $0.9<y<0.95$  for injection with $q = 10^{-3}$. 
\end{itemize}
The total simulation time is $t_{max} = 0.1$, and is discretized into 20 time steps for both upscaled and fine-scale solvers.

\begin{figure}[h!]
\centering
\includegraphics[width=1 \textwidth]{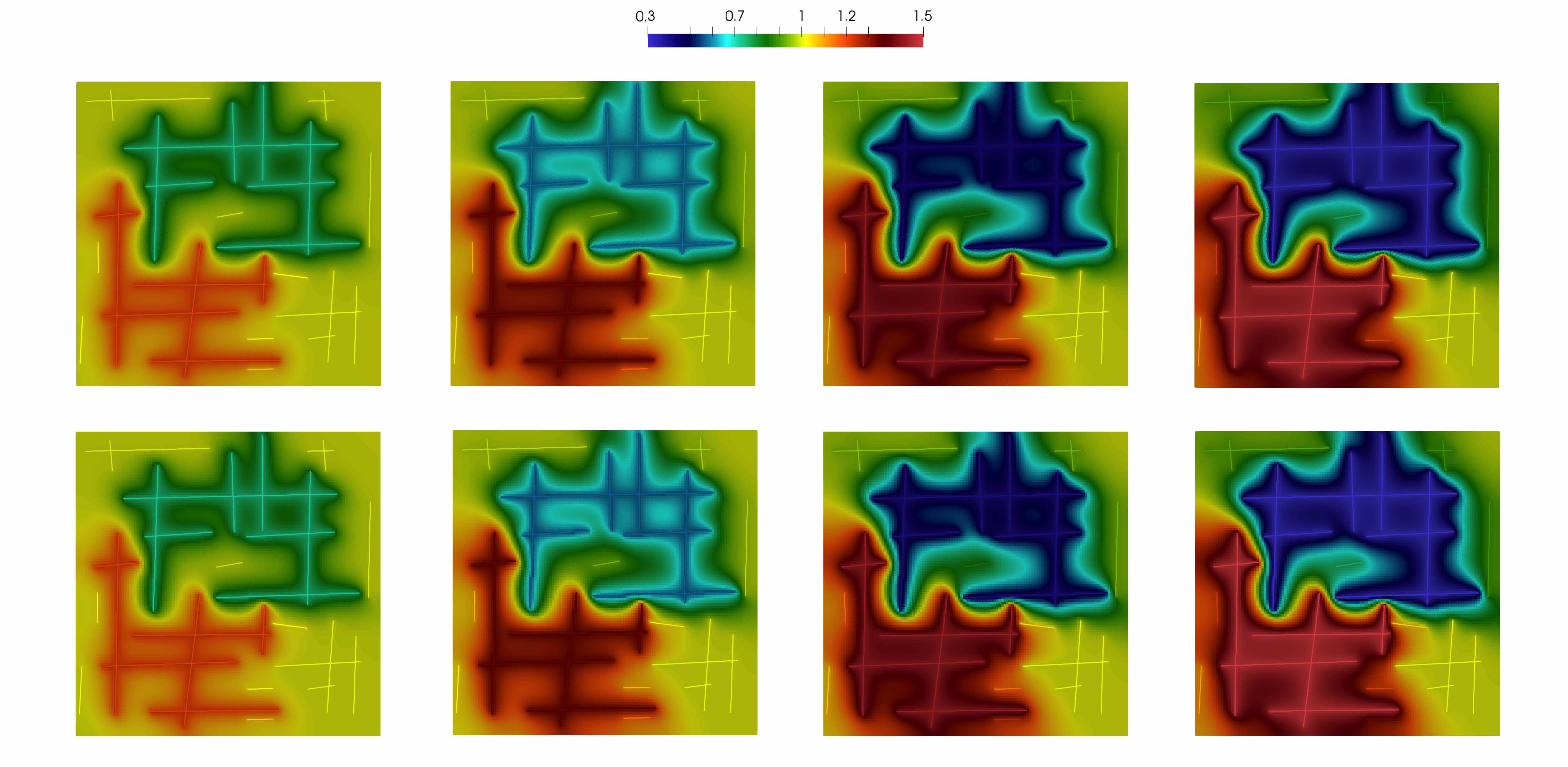}
\caption{EFM and DFM fine-scale solutions. Geometry 1 for different time steps $t_5 = 0.025$ (first column), $t_{10} = 0.05$ (second column), $t_{15} = 0.075$ (third column) and $t_{20} = 0.1$ (fourth column). First row: DFM fine-scale solution. Second row: EFM fine-scale solution}
\label{fig:uf-t1}
\end{figure}

\begin{figure}[h!]
\centering
\includegraphics[width=1 \textwidth]{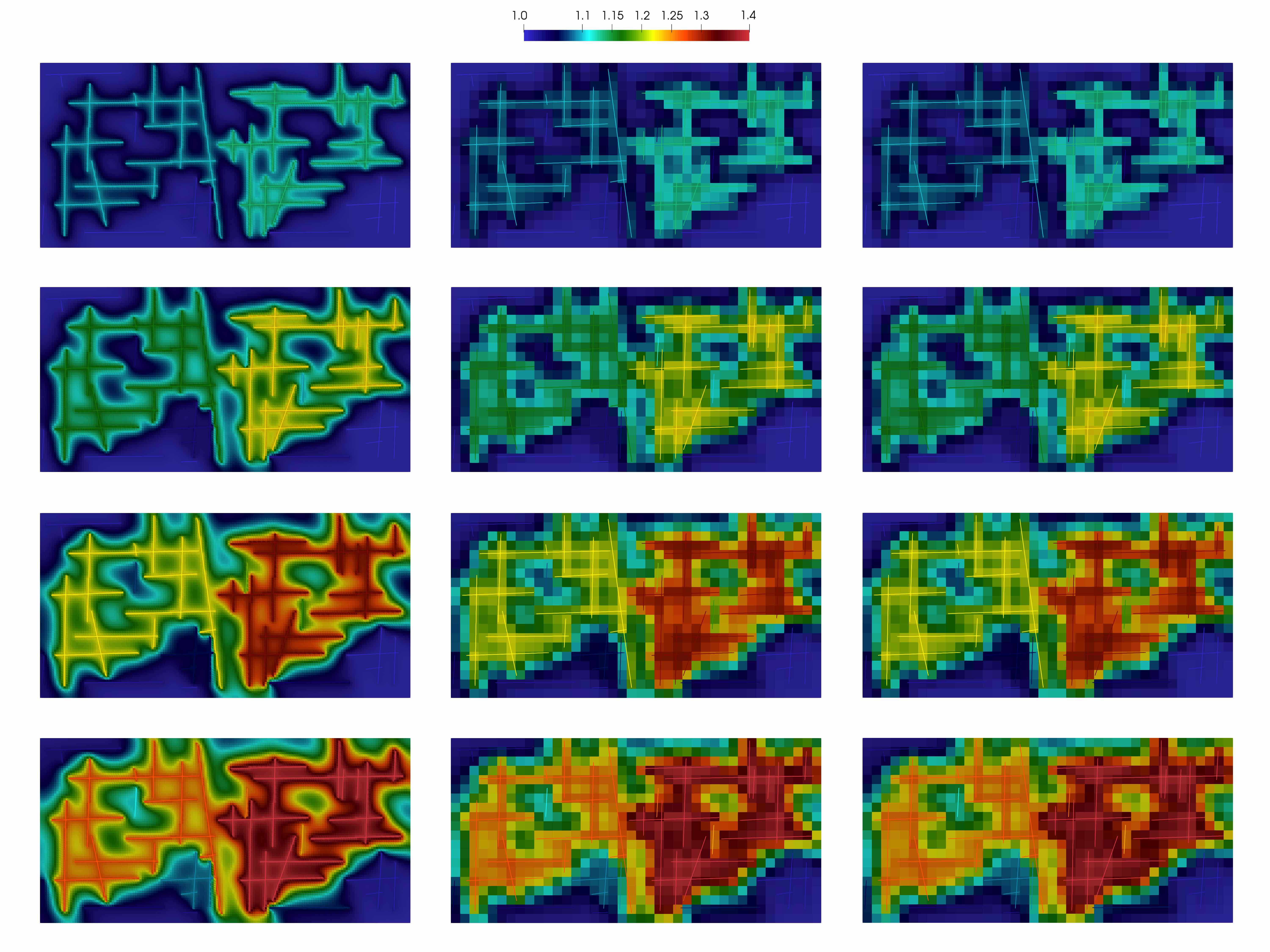}
\caption{Multiscale solutions on mesh $40 \times 20$ with $K^2$ using EFM fine-scale solver. Geometry 2 for different time steps $t_5 = 0.025$, $t_{10} = 0.05$, $t_{15} = 0.075$ and $t_{20} = 0.1$ (from top to bottom). First column: fine scale solution. Second column: cell average  for fine grid solution. Third column: upscaled coarse grid solution}
\label{fig:u-efm-t2}
\end{figure}

\begin{table}[h!]
\centering
\begin{tabular}{ |c | c | c | c | c | }
\hline
$K^s$ & 
$t_5$ & 
$t_{10}$ & 
$t_{15}$ & 
$t_{20}$ \\ \hline
\multicolumn{5}{|c|}{Geometry 1 with DFM. $20 \times 20$} \\
\hline
$s = 1$	& 0.433	&	0.955	&	1.623	&	2.404 \\ \hline
$s = 2$	& 0.233	&	0.314	&	0.396	&	0.476 \\ \hline
$s = 3$	& 0.201	&	0.305	&	0.378	&	0.455 \\ \hline
\multicolumn{5}{|c|}{Geometry 2 with DFM. $40 \times 20$} \\
\hline
$s = 2$	& 0.138	&	0.170	&	0.192	&	0.211 \\ \hline
\end{tabular}
\,\,\,\,
\begin{tabular}{ |c | c | c | c | c | }
\hline
$K^s$ & 
$t_5$ & 
$t_{10}$ & 
$t_{15}$ & 
$t_{20}$ \\ \hline
\multicolumn{5}{|c|}{Geometry 1 with EFM. $20 \times 20$} \\
\hline
$s = 1$	& 0.380	&	0.724	&	1.127	&	1.578 \\ \hline
$s = 2$	& 0.264	&	0.390	&	0.499	&	0.598 \\ \hline
$s = 3$	& 0.259	&	0.388	&	0.483	&	0.579 \\ \hline
\multicolumn{5}{|c|}{Geometry 2 with EFM. $40 \times 20$} \\
\hline
$s = 2$	& 0.098	&	0.133	&	0.158	&	0.180 \\ \hline
\end{tabular}
\caption{Relative errors of the average cell solution on a coarse mesh. Right: DFM. Left: EFM. Geometry 1 and 2. }
\label{err-t12}
\end{table}

We first consider Geometry 1. For DFM model, the unstructured fine grid contains 47520 fine-scale elements for porous matrix and 1042 fine-scale elements for fractures. For EFM model, we use the structured fine grid containing 20000 fine-scale elements (matrix) and 1042 fine-scale elements (fractures). We consider uniformly structured coarse grid with 400 coarse-scale elements($20 \times 20$).
Fine-grid solutions using DFM and EFM models are presented in Figure \ref{fig:uf-t1}. We notice similar solutions for both models for sufficient fine grids.  

In Table \ref{err-t12}, we show relative errors for different number of oversampling layers $K^s$ with $s = 1,2$ and $3$, using DFM and EFM fine grid approximations.   
For coarse mesh with 400 cells, when we take 2 oversampling layers, we have $0.398 \%$ of error at final time. 
The fine-scale systems have $DOF_f = 47520 + 1042$ for DFM and $DOF_f = 20000 + 1042$ for EFM. While the upscaled model only has $DOF_c = 593$ for coarse mesh with 400 cells ($20 \times 20$).
We note that our proposed method provide with huge reduction of the system size and very accurate approximations. 

Next, we consider Geometry 2. For DFM model, we use unstructured fine grid containing 98398 fine cells (matrix) and 2170 fine cells (fractures). For EFM mode, we employ structured fine grid containing 25600 fine cells (matrix) and 2170 fine cells (fractures).
In this test, we consider uniform structured coarse grids $40 \times 20$ (800 cells).
The fine-scale systems have $DOF_f = 98398 + 2170$ for DFM and $DOF_f = 25600 + 2170$ for EFM. The upscaled model  has $DOF_c = 1179$ for coarse mesh with 800 cells ($40 \times 20$).
In Figure \ref{fig:u-efm-t2}, we present results using the upscaled model, where we obtain very accurate simulation results with very small DOF in the upscaled model. For local domain $K^2$, relative error for porous matrix are less than one percent.

In Table \ref{err-t12}, we present relative errors between the solutions of EFM and DFM fine-grid models and the upscaled model, where we only use two oversampling layers $K^2$ for basis construction.
From the numerical results, we observe a good convergence when we take sufficient number of oversampled layers.

\section{Numerical results for heterogeneous permeability}
\label{sec:num3}

\begin{figure}[h!]
\centering
\includegraphics[width=0.5 \textwidth]{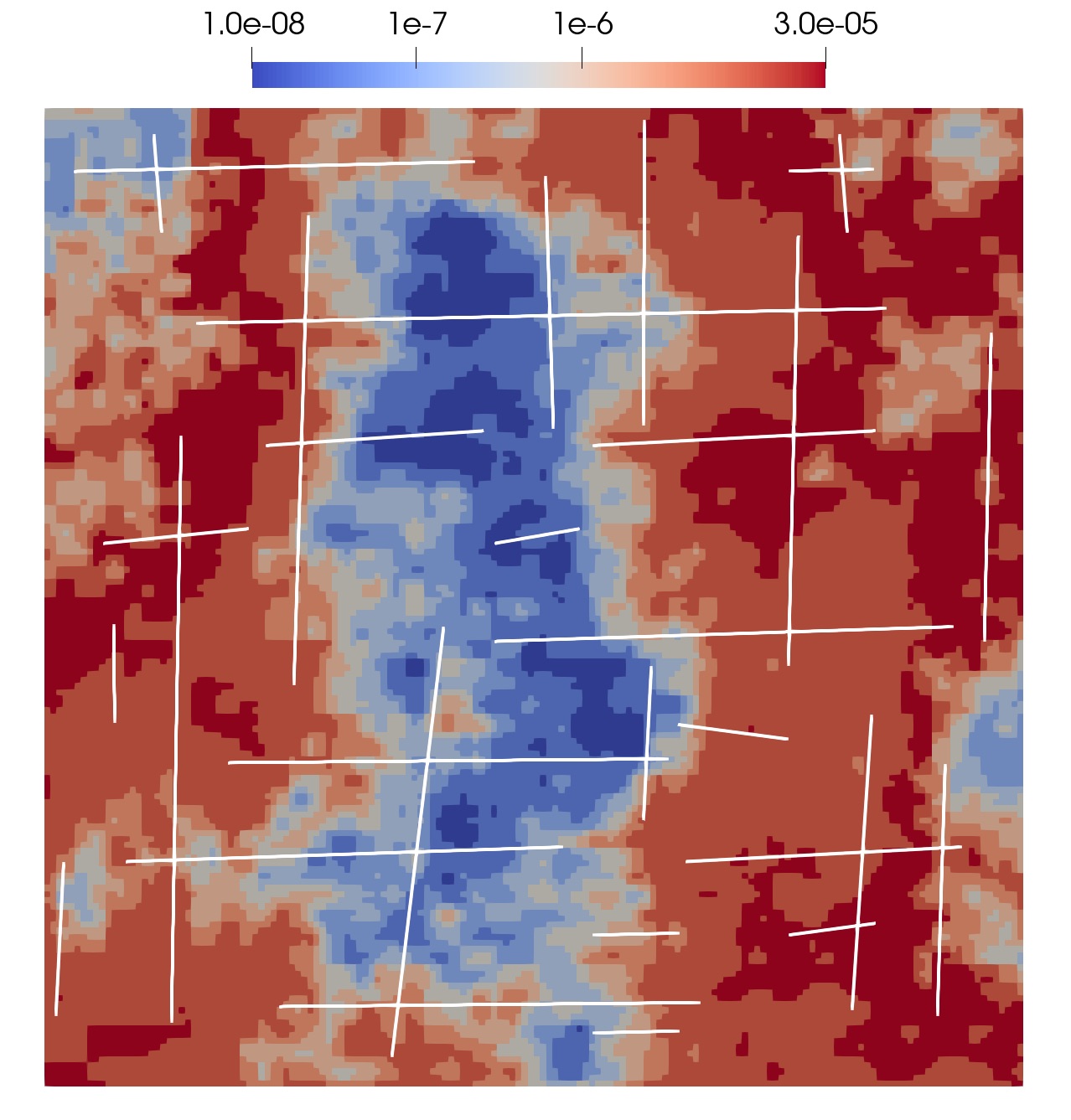}
\caption{Heterogeneous pore matrix permeability for Geometry 1} 
\label{fig:km}
\end{figure}

Finally, we consider a test case with heterogeneous permeability (Figure \ref{fig:km}) for Geometry 1 with same parameters as in the previous section. As for fine grid approximation, we use EFM. The simulation time is $t_{max} = 0.1$ with 100 time steps. 

\begin{figure}[h!]
\centering
\includegraphics[width=1 \textwidth]{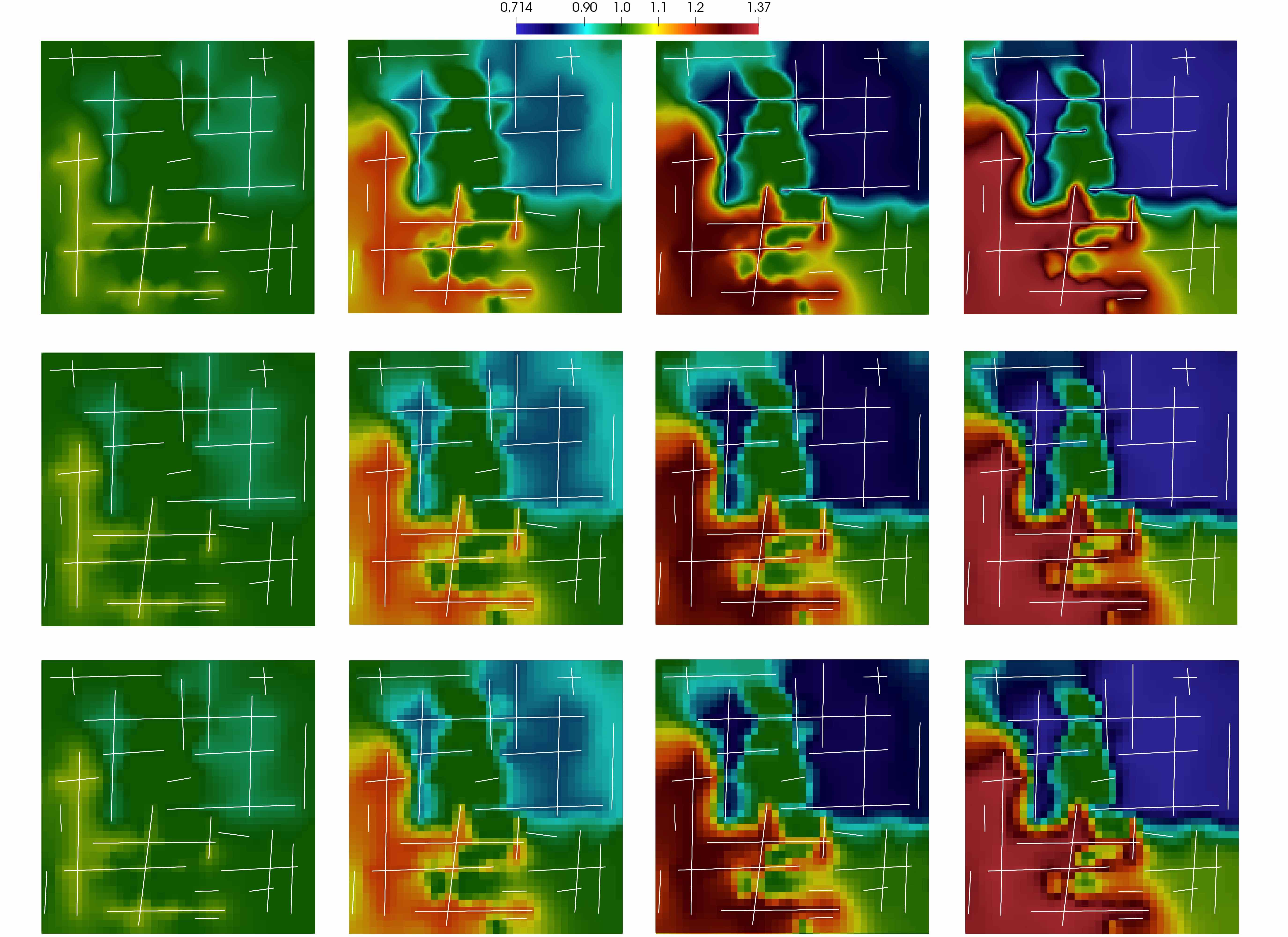}
\caption{Multiscale solutions on mesh $40 \times 40$ with $K^2$ using EFM fine-scale solver. Geometry 1 for different time steps $t_5$, $t_{25}$, $t_{50}$ and $t_{100}$ (from left to right). 
First row: fine scale solution. Second row: cell average  for fine grid solution. Third row: upscaled coarse grid solution}
\label{fig:u-km}
\end{figure}

\begin{table}[h!]
\centering
\begin{tabular}{ |c | c | c | c | c | }
\hline
$K^s$ & 
$t_5$ & 
$t_{25}$ & 
$t_{50}$ & 
$t_{100}$ \\ \hline
\multicolumn{5}{|c|}{Geometry 1. $40 \times 40$} \\
\hline
$s = 1$	& 0.137 &	0.434	&	0.581	&	0.785 \\ \hline
$s = 2$	& 0.078 &	0.182  &	0.230  &	0.306 \\ \hline
$s = 3$	& 0.078  &	0.178  &	0.200  &	0.186 \\ \hline
$s = 4$	& 0.078  &	0.176  &	0.198  &	0.183 \\ \hline
\end{tabular}
\caption{Relative errors of the average cell solution on a coarse mesh. Geometry 1 with heterogeneous permeability. }
\label{err-km}
\end{table}

In this test, we consider uniform structured coarse grids $40 \times 40$ (1600 cells).
In Figure \ref{fig:u-km}, we present results using the upscaled model, where relative error are less than one percent.
In Table \ref{err-km}, we present relative errors between the fine grid and the upscaled coarse grid solutions. 
The fine-scale systems have size $DOF_f = 27384$ with solution time $57$ seconds. The upscaled model has $DOF_c = 1165$ with solution time $5.9$ seconds. The computational time is reduced 10 times due to the reduction in the size of the system. The proposed method is shown to be very efficient and provides good accuracy.

\section{Conclusion} \label{sec:conclusion}

We consider mixed dimensional coupled problem for flow simulation in fractured porous media for EFM and DFM with finite volume approximation on the fine grid.
We presented an upscaling method for coupled problems in fractured domains. In this work, we construct multiscale basis function for background medium and additional multiscale basis for fractures.

We presented numerical results for model problems: (1) discrete fracture fine grid model with low and high permeable fractures; (2) embedded fine grid model for two geometries with different number of fracture lines and (3) embedded fracture fine grid model  with heterogeneous permeability. 
Our proposed upscaling method provides good accuracy and give a significant reduction in the size of  the problem system. The resulting upscaled model has minimal size and the solution obtained has physical meaning on the coarse grid.

\section{Acknowledgements}
MV's  work is supported by the grant of the Russian Scientific Found N17-71-20055.
VA's work is supported by the mega-grant of the Russian Federation Government (N 14.Y26.31.0013). 
EC's work is partially supported by Hong Kong RGC General Research Fund (Project 14304217)
and CUHK Direct Grant for Research 2017-18.

\bibliographystyle{plain}
\bibliography{lit}

\begin{thebibliography}{10}

\bibitem{akkutlu2018multiscale}
I~Yucel Akkutlu, Yalchin Efendiev, Maria Vasilyeva, and Yuhe Wang.
\newblock Multiscale model reduction for shale gas transport in poroelastic
  fractured media.
\newblock {\em Journal of Computational Physics}, 353:356--376, 2018.

\bibitem{akkutlu2015multiscale}
IY~Akkutlu, Yalchin Efendiev, and Maria Vasilyeva.
\newblock Multiscale model reduction for shale gas transport in fractured
  media.
\newblock {\em Computational Geosciences}, pages 1--21, 2015.

\bibitem{barenblatt1960basic}
GI~Barenblatt, Iu~P Zheltov, and IN~Kochina.
\newblock Basic concepts in the theory of seepage of homogeneous liquids in
  fissured rocks [strata].
\newblock {\em Journal of applied mathematics and mechanics}, 24(5):1286--1303,
  1960.

\bibitem{bosma2017multiscale}
Sebastian Bosma, Hadi Hajibeygi, Matei Tene, and Hamdi~A Tchelepi.
\newblock Multiscale finite volume method for discrete fracture modeling on
  unstructured grids (ms-dfm).
\newblock {\em Journal of Computational Physics}, 2017.

\bibitem{CELV2015}
E.~T. Chung, Y.~Efendiev, G.~Li, and M.~Vasilyeva.
\newblock Generalized multiscale finite element method for problems in
  perforated heterogeneous domains.
\newblock {\em to appear in Applicable Analysis}, 255:1--15, 2015.

\bibitem{chung2016adaptive}
Eric Chung, Yalchin Efendiev, and Thomas~Y Hou.
\newblock Adaptive multiscale model reduction with generalized multiscale
  finite element methods.
\newblock {\em Journal of Computational Physics}, 320:69--95, 2016.

\bibitem{chung2017coupling}
Eric~T Chung, Yalchin Efendiev, Tat Leung, and Maria Vasilyeva.
\newblock Coupling of multiscale and multi-continuum approaches.
\newblock {\em GEM-International Journal on Geomathematics}, 8(1):9--41, 2017.

\bibitem{chung2017constraint}
Eric~T Chung, Yalchin Efendiev, and Wing~Tat Leung.
\newblock Constraint energy minimizing generalized multiscale finite element
  method.
\newblock {\em arXiv preprint arXiv:1704.03193}, 2017.

\bibitem{chung2017non}
Eric~T Chung, Yalchin Efendiev, Wing~Tat Leung, Yating Wang, and Maria
  Vasilyeva.
\newblock Non-local multi-continua upscaling for flows in heterogeneous
  fractured media.
\newblock {\em arXiv preprint arXiv:1708.08379}, 2017.

\bibitem{Quarteroni2008coupling}
Carlo D'angelo and Alfio Quarteroni.
\newblock On the coupling of 1d and 3d diffusion-reaction equations:
  application to tissue perfusion problems.
\newblock {\em Mathematical Models and Methods in Applied Sciences},
  18(08):1481--1504, 2008.

\bibitem{douglas1990dual}
Jim Douglas~Jr and T~Arbogast.
\newblock Dual porosity models for flow in naturally fractured reservoirs.
\newblock {\em Dynamics of Fluids in Hierarchical Porous Media}, pages
  177--221, 1990.

\bibitem{d2012mixed}
Carlo D’Angelo and Anna Scotti.
\newblock A mixed finite element method for darcy flow in fractured porous
  media with non-matching grids.
\newblock {\em ESAIM: Mathematical Modelling and Numerical Analysis},
  46(2):465--489, 2012.

\bibitem{weinan2007heterogeneous}
Weinan E, Bjorn Engquist, Xiantao Li, Weiqing Ren, and Eric Vanden-Eijnden.
\newblock Heterogeneous multiscale methods: a review.
\newblock {\em Commun. Comput. Phys}, 2(3):367--450, 2007.

\bibitem{EGG_MultiscaleMOR}
Y.~Efendiev, J.~Galvis, and E.~Gildin.
\newblock Local-global multiscale model reduction for flows in highly
  heterogeneous media.
\newblock {\em Journal of Computational Physivs}, 231 (24):8100--8113, 2012.

\bibitem{egh12}
Y.~Efendiev, J.~Galvis, and T.~Hou.
\newblock Generalized multiscale finite element methods.
\newblock {\em Journal of Computational Physics}, 251:116--135, 2013.

\bibitem{eh09}
Y.~Efendiev and T.~Hou.
\newblock {\em {Multiscale Finite Element Methods: Theory and Applications}},
  volume~4 of {\em Surveys and Tutorials in the Applied Mathematical Sciences}.
\newblock Springer, New York, 2009.

\bibitem{efendiev2015hierarchical}
Yalchin Efendiev, Seong Lee, Guanglian Li, Jun Yao, and Na~Zhang.
\newblock Hierarchical multiscale modeling for flows in fractured media using
  generalized multiscale finite element method.
\newblock {\em arXiv preprint arXiv:1502.03828}, 2015.
\newblock to appear in International Journal on Geomathematics, (DOI)
  10.1007/s13137-015-0075-7.

\bibitem{formaggia2014reduced}
Luca Formaggia, Alessio Fumagalli, Anna Scotti, and Paolo Ruffo.
\newblock A reduced model for darcy’s problem in networks of fractures.
\newblock {\em ESAIM: Mathematical Modelling and Numerical Analysis},
  48(4):1089--1116, 2014.

\bibitem{garipov2016discrete}
TT~Garipov, M~Karimi-Fard, and HA~Tchelepi.
\newblock Discrete fracture model for coupled flow and geomechanics.
\newblock {\em Computational Geosciences}, 20(1):149--160, 2016.

\bibitem{ginting2011application}
Victor Ginting, Felipe Pereira, Michael Presho, and Shaochang Wo.
\newblock Application of the two-stage markov chain monte carlo method for
  characterization of fractured reservoirs using a surrogate flow model.
\newblock {\em Computational Geosciences}, 15(4):691, 2011.

\bibitem{hkj12}
H.~Hajibeygi, D.~Kavounis, and P.~Jenny.
\newblock A hierarchical fracture model for the iterative multiscale finite
  volume method.
\newblock {\em Journal of Computational Physics}, 230(24):8729--8743, 2011.

\bibitem{hoteit2008efficient}
Hussein Hoteit and Abbas Firoozabadi.
\newblock An efficient numerical model for incompressible two-phase flow in
  fractured media.
\newblock {\em Advances in Water Resources}, 31(6):891--905, 2008.

\bibitem{houwu97}
T.~Hou and X.H. Wu.
\newblock A multiscale finite element method for elliptic problems in composite
  materials and porous media.
\newblock {\em J. Comput. Phys.}, 134:169--189, 1997.

\bibitem{jenny2005adaptive}
Patrick Jenny, Seong~H Lee, and Hamdi~A Tchelepi.
\newblock Adaptive multiscale finite-volume method for multiphase flow and
  transport in porous media.
\newblock {\em Multiscale Modeling \& Simulation}, 3(1):50--64, 2005.

\bibitem{karimi2003efficient}
Mohammad Karimi-Fard, Luis~J Durlofsky, Khalid Aziz, et~al.
\newblock An efficient discrete fracture model applicable for general purpose
  reservoir simulators.
\newblock In {\em SPE Reservoir Simulation Symposium}. Society of Petroleum
  Engineers, 2003.

\bibitem{karimi2001numerical}
Mohammad Karimi-Fard, Abbas Firoozabadi, et~al.
\newblock Numerical simulation of water injection in 2d fractured media using
  discrete-fracture model.
\newblock In {\em SPE annual technical conference and exhibition}. Society of
  Petroleum Engineers, 2001.

\bibitem{logg2009efficient}
Anders Logg.
\newblock Efficient representation of computational meshes.
\newblock {\em International Journal of Computational Science and Engineering},
  4(4):283--295, 2009.

\bibitem{logg2012automated}
Anders Logg, Kent-Andre Mardal, and Garth Wells.
\newblock {\em Automated solution of differential equations by the finite
  element method: The FEniCS book}, volume~84.
\newblock Springer Science \& Business Media, 2012.

\bibitem{lunati2006multiscale}
Ivan Lunati and Patrick Jenny.
\newblock Multiscale finite-volume method for compressible multiphase flow in
  porous media.
\newblock {\em Journal of Computational Physics}, 216(2):616--636, 2006.

\bibitem{martin2005modeling}
Vincent Martin, J{\'e}r{\^o}me Jaffr{\'e}, and Jean~E Roberts.
\newblock Modeling fractures and barriers as interfaces for flow in porous
  media.
\newblock {\em SIAM Journal on Scientific Computing}, 26(5):1667--1691, 2005.

\bibitem{schwenck2015dimensionally}
Nicolas Schwenck, Bernd Flemisch, Rainer Helmig, and Barbara~I Wohlmuth.
\newblock Dimensionally reduced flow models in fractured porous media:
  crossings and boundaries.
\newblock {\em Computational Geosciences}, 19(6):1219--1230, 2015.

\bibitem{tene2016multiscale}
M~Tene, MS~Al~Kobaisi, and H~Hajibeygi.
\newblock Multiscale projection-based embedded discrete fracture modeling
  approach (f-ams-pedfm).
\newblock In {\em ECMOR XV-15th European Conference on the Mathematics of Oil
  Recovery}, 2016.

\bibitem{ctene2016algebraic}
Matei {\c{T}}ene, Mohammed~Saad Al~Kobaisi, and Hadi Hajibeygi.
\newblock Algebraic multiscale method for flow in heterogeneous porous media
  with embedded discrete fractures (f-ams).
\newblock {\em Journal of Computational Physics}, 321:819--845, 2016.

\bibitem{ctene2017projection}
Matei {\c{T}}ene, Sebastian~BM Bosma, Mohammed~Saad Al~Kobaisi, and Hadi
  Hajibeygi.
\newblock Projection-based embedded discrete fracture model (pedfm).
\newblock {\em Advances in Water Resources}, 105:205--216, 2017.

\bibitem{warren1963behavior}
JE~Warren, P~Jj Root, et~al.
\newblock The behavior of naturally fractured reservoirs.
\newblock {\em Society of Petroleum Engineers Journal}, 3(03):245--255, 1963.

\end{thebibliography}

\end{document}